\documentclass[letterpaper,10pt,twocolumn]{article} 
\usepackage{amsmath}
\usepackage{amssymb} 
\usepackage{mathrsfs}
\usepackage{verbatim}
\usepackage{graphicx}
\usepackage{algorithm,algorithmic}
\usepackage{enumitem}
\usepackage{cite}
\usepackage[margin=1in]{geometry}

\newtheorem{assumption}{Assumption}
\newtheorem{theorem}{Theorem}
\newtheorem{lemma}[theorem]{Lemma}
\newtheorem{prop}[theorem]{Proposition}

\newtheorem{definition}{Definition}
\newtheorem{problem}{Problem}
\newtheorem{example}{Example}
\newtheorem{remark}{Remark}

\newcommand{\NN}{\mathbb{N}}
\newcommand{\RR}{\mathbb{R}}
\newcommand{\CC}{\mathbb{C}}

\newcommand{\cC}{\mathcal{C}}

\newcommand{\cL}{\mathcal{L}}
\newcommand{\cLinf}{\mathcal{L}_{\infty}}

\newcommand{\cP}{\mathcal{P}}

\newcommand{\cS}{\mathcal{S}}

\newcommand{\cV}{\mathcal{V}}
\newcommand{\cW}{\mathcal{W}}
\newcommand{\cWp}{\mathcal{W}^{\perp}}

\newcommand{\al}{\alpha}

\newcommand{\ale}{{\alpha}_e}
\newcommand{\be}{{\beta}_e}
\newcommand{\ct}{\tilde{c}_i}

\newcommand{\eps}{\epsilon}
\newcommand{\stol}{\sigma_{\text{tol}}}

\newcommand{\shattol}{\hat{\sigma}_{\text{tol}}}
\newcommand{\veps}{\varepsilon}

\newcommand{\hci}{\hat{c}_{i}}

\newcommand{\dhci}{\dot{\hat{c}}_{i}}

\newcommand{\Fb}{\bar{F}}
\newcommand{\Fbz}{\bar{F_{\circ}}}

\newcommand{\Gb}{\bar{G}}

\newcommand{\nh}{\hat{n}}
\newcommand{\nb}{\bar{n}}

\newcommand{\psih}{\hat{\psi}}
\newcommand{\psit}{\tilde{\psi}}
\newcommand{\psis}{\psi^{\star}}
\newcommand{\rnom}{r^{\circ}_i}
\newcommand{\rhat}{\hat{r}_i}

\newcommand{\dsigr}{\dot{\hat{\sigma}}_{i}}
\newcommand{\sigr}{\hat{\sigma}_{i}}
\newcommand{\sigt}{\tilde{\sigma}_{i}}

\newcommand{\ths}{\theta^{\star}}

\newcommand{\tho}{\theta_{\circ}}

\newcommand{\wb}{\bar{w}}

\newcommand{\xihb}{\hat{\xi}_{i+1:}}
\newcommand{\xihbi}{\hat{\xi}_{i:}}
\newcommand{\xibt}{\tilde{\xi}_{i+1:}}
\newcommand{\xih}{\hat{\xi}}

\newcommand{\xis}{{\xi}^{\star}} 
\newcommand{\xisb}{\xi^{\star}_{i+1:}}
\newcommand{\xisbi}{\xi^{\star}_{i:}}
\newcommand{\xit}{\tilde{\xi}}
\newcommand{\yb}{\bar{y}}

\newcommand{\pe}{_{pe}}

\newcommand{\inv}{^{-1}}
\newcommand{\0}{_{\circ}}

\newcommand{\T}{^{\intercal}}

\providecommand{\Proj}{\mathrm{Proj}} 
 
\newcommand{\Zero}{\mathbf{0}}

\newcommand{\AVG}[1]{\lim_{T \to \infty} \frac{1}{T} \int_{t_0}^{t_0 + T} #1 \,d\tau}

\newcommand{\tqed}{\hfill$\triangleleft$}

\DeclareMathOperator{\Img}{Im}

\DeclareMathOperator{\Ker}{Ker}

\DeclareMathOperator{\diag}{diag}

\DeclareMathOperator{\rnk}{rank}

\DeclareMathOperator{\interior}{int}
\DeclareMathOperator{\eig}{eig}

\usepackage{subcaption}
\usepackage{graphicx}
\usepackage{url}

\begin{document}

\title{Adaptive Identification via Regressor Redundancy: Extended Version}
\author{
Fatima J. Ghadieh$^{1}$ \and
Mireille E. Broucke$^{1,2}$
\thanks{$^{1}$Supported by the Natural Sciences and Engineering Research Council of Canada (NSERC). Electrical and Computer Engineering, University of Toronto, Toronto, ON, Canada.}%
\thanks{$^{2}$Email: \texttt{broucke@control.utoronto.ca}}%
}
\maketitle

\begin{abstract} 
We develop a continuous time adaptive identification scheme for linear systems of uncertain order. A key enabling concept is that of regressor redundancy, which tests whether a component of an adaptive observer state adds new excitation to the remaining components. We establish that overparameterization of the observer causes a structured loss of persistent excitation, and we exploit this structure to identify in real time the model order simultaneously with the plant parameters.
\end{abstract}
 
\section{Introduction} 

In this paper we continue our investigation of system identification methods that operate purely in continuous time \cite{GHADIEH26}, motivated by studies in systems neuroscience where we wish to understand how the brain learns dynamic models \cite{BROUCKE26}.
Adaptive observers offer a natural framework for online identification in continuous time \cite{KREISSEL77}, but they require the plant order to be known \cite{LAI24}. When only an upper bound on the order is known, the associated regressor is no longer persistently exciting (PE) \cite{KREISSEL96}, and the plant parameters cannot be uniquely recovered by a standard adaptive law. This raises the question of whether the loss of PE caused by overparameterization has a structure that can be exploited to identify the plant order.

The model order is usually estimated offline from batch sampled data, either by minimizing an information criterion \cite{LJUNG99} or by determining the rank of a (Hankel) matrix constructed from input-output data, typically through a singular value decomposition (SVD) \cite{HOKALMAN66,VANOVERSCHEE96}; see also \cite{SARKAR21,ZHENG25} for recent work. Online order estimation has mostly been studied in discrete time. In \cite{BAKO09} the order is recursively tracked by testing the singularity of a data covariance matrix without using an SVD. The paper \cite{HASHAMBHOY05} runs two identifiers with different initial conditions in order to recursively reduce the order bound. The continuous time case is less developed; in \cite{HUANG22} the order is estimated online from sampled data and with an SVD at each update.  Our previous work \cite{GHADIEH26} determined the order with an incremental continuous-time Ho-Kalman-like algorithm based on an SVD, with parameters computed in a separate step.  In contrast, the present paper proposes a design purely in continuous time that requires neither sampled data nor an SVD, and it identifies the order and the parameters simultaneously in a single run. 

We show using the adaptive observer of \cite{KREISSEL77} that
overparameterization causes a structured loss of PE in the regressor that reveals the model order. We introduce the notion of \emph{regressor redundancy} along with a test that detects novel excitation in regressor components. The test informs a parameter adaptation law about whether each component of the adapted parameter should be learned or forgotten.  A related mechanism appears in \cite{MARINO07}, where the goal is output regulation with an unknown exosystem. There, the order of the exosystem is identified by detecting a loss of PE through determinants computed from a bank of filters. In contrast,  our scheme is computationally attractive: it uses only linear filters and gradient laws, and it requires no determinants, matrix decompositions such as the SVD, inversions, or sampling of input-output data.

\section{Background}
We use the following standard notation. 
The vector with $p$ zeros is denoted $\Zero_p$.
The disjoint union of two sets is denoted $\sqcup$. 
For a subspace $\cV \subset \RR^n$, its dimension is denoted $\dim(\cV)$. 
For a matrix $A$, its image and kernel are denoted $\Img(A)$ and $\Ker(A)$. For a square matrix $A \in \RR^{n \times n}$, its spectrum, rank, and transpose are denoted as $\eig(A)$, $\rnk(A)$, and $A\T$. The smallest and largest eigenvalues of a symmetric matrix $A \in \RR^{n \times n}$ are denoted as $\lambda_{\min}(A)$ and $\lambda_{\max}(A)$.
The degree of polynomial $f$ is denoted as $\deg(f)$. 
Let $A\0$ be the matrix of all zero elements, except for ones on the upper diagonal, and let $B\0 = (0,\ldots,0,1)$ (dimensions are inferred from the context). If $u(t)$ is a function of time, then $u(s)$ denotes its Laplace transform. For a time-varying signal $w: [0, \infty) \to \RR^n$, we say $w \in \cL_2$ if it is square-integrable, i.e. $\int_0^\infty \|w(t)\|^2 dt < \infty$. We say $w \in \cLinf$ if it is uniformly bounded, i.e. $\sup_{t \ge 0} \|w(t)\| < \infty$. Moreover, let $\| w \|_{\cL_\infty}$ denote the $\cL_\infty$ norm of a signal $w(t) \in \RR^n$ (assuming it exists). 
A bounded signal $w(t) \in \RR^q$ is {\em persistently exciting} (PE) if there exist $\beta\pe, T\pe > 0$ such that
\begin{equation}\label{eq:PE}
\frac{1}{T\pe} \int_{t}^{t+T\pe} w(\tau) w\T(\tau) \,d\tau  
\succeq \beta\pe I \,, \quad \forall \ t \geq t_0 \,.
\end{equation}
We call the pair $(\beta\pe,T\pe)$ {\em PE constants of $w$} (they are not unique). 

\begin{definition} \label{def:auto}
A signal $w(t) \in \RR^q$ is said to be {\em stationary} if it is bounded and piecewise continuous, and its {\em autocovariance matrix}
$R_w( \rho ) := \AVG{w(\tau) w\T(\tau + \rho )}$
exists, is independent of the initial time $t_0$, and the convergence of the time average is uniform in $t_0 \geq 0$. 
\end{definition} 
We only require the zero lag autocovariance, denoted as $R[w] = R_w(0)$. If a signal $w:\RR \rightarrow \RR$ has a sinusoidal component at frequency $\omega_0$, then its {\em spectral measure} $S_w(\omega)$ has \textit{point masses} at $\omega_0$ and $-\omega_0$. A stationary signal $w(t)$ is \textit{sufficiently rich} (SR) of order $k$ if the support of $S_w(\omega)$ contains at least $k$ points \cite[p.~255]{IOANNOU12}.

\subsection{PE Subspaces}\label{sec:PEsubspace}

Consider a stationary signal $w(t) \in \RR^q$. The {\em PE subspace} of $w$ is defined to be $\cW = \Img ( R[w] )$, while the {\em non-PE subspace} of $w$ is $\cW^{\perp}$.  We denote $q\pe := \dim(\cW)$ as its {\em degree of persistent excitation}.  When $q\pe = q$, $w$ is PE by \cite[Proposition~2.7.1]{SASTRY89}, whereas when $q\pe = 0$, $w$ has {\em no persistent excitation}.  

\begin{lemma}\label{lem:kerR}
Let $w(t) \in \RR^{q}$ be a stationary signal with autocovariance $R[w]$ and PE subspace $\cW = \Img(R[w])$. For every $a \in \RR^{q}$, $a\T w$ is not PE if and only if $a \in \cWp$. 
\tqed
\end{lemma}

\section{Problem Statement}

Consider a SISO LTI system modeled as a transfer function
\begin{align}
\label{eq:theplant}
H(s) &= \frac{\beta(s)}{\al(s)} =: 
\frac{\beta_{n-1} s^{n-1} + \ldots + \beta_1 s + \beta_0} 
     {s^n + \al_{n-1} s^{n-1}  + \ldots + \al_1 s + \al_0} \,,
\end{align}
where $n$ is the order of the plant. 
\begin{assumption}
\label{assum:plant}
The system \eqref{eq:theplant} and the input $u$ satisfy the following.
\begin{enumerate}[label=(A\arabic*), leftmargin=*, ref=(A\arabic*)]
\item \label{as:coprime}
$\{ \al(s), \beta(s) \}$ are coprime.
\item \label{as:hurwitz}
 $\al(s)$ is Hurwitz.
\item \label{as:unknown}
The model order $n$ and parameters $\al := ( \al_0, \ldots, \al_{n-1} ) \in \RR^n$ and $\beta := ( \beta_0, \ldots, \beta_{n-1} ) \in \RR^n$ are unknown.
\item \label{assum:over} 
An upper bound on the plant order $\nb \ge n$ is known.
Let $p := \nb - n$ be the excess order. 
\item \label{as:measurements}
The measurements are input $u(t) \in \RR$ and output $y(t) \in \RR$.
\item \label{assum:input}
The input $u(t)$ is a stationary multisine signal that is sufficiently rich of order $2\nb$. 
\tqed
\end{enumerate}
\end{assumption}

\begin{problem}
Consider the system \eqref{eq:theplant} and suppose Assumption~\ref{assum:plant} holds. Develop a continuous time procedure to simultaneously identify the plant order $n$ in finite time and parameters $(\al,\beta)$ exponentially.
\tqed
\end{problem} 

\section{Overparameterized Model}

Consider the {\em overparameterized model} 
\begin{subequations}
\label{eq:overp1}
\begin{align}
\dot{\xi}_1 &= F \xi_1 + G y\0 \\
\dot{\xi}_2 &= F \xi_2 + G u \\
y\0 & = \theta_1\T \xi_1 + \theta_2\T \xi_2 \,, 
\end{align}
\end{subequations}
where $\xi := (\xi_1, \xi_2)$ with $\xi_i \in \RR^{\nb}$ is the state, $y\0 \in \RR$ is the overparameterized model output, and $\theta_1, \theta_2 \in \RR^{\nb}$ are to-be-defined parameters. We assume $F$ is Hurwitz, and $(F,G)$ is a controllable pair. Specifically, we select
\begin{align}
\label{eq:FG}
F &= A\0 + B\0 (-f)\T \,, \qquad G = B\0 \,, 
\end{align}
where $f = (f_0, f_1, \ldots, f_{\nb-1}) \in \RR^{\nb}$ and 
$f(s) = s^{\nb} + f_{\nb-1} s^{\nb-1} + \cdots + f_1 s + f_0$ is a Hurwitz polynomial. Defining $\theta := (\theta_1,\theta_2) \in \RR^{2 \nb}$, \eqref{eq:overp1} can be written compactly as 
\begin{subequations}
\label{eq:overp2}
\begin{align}
\nonumber
\dot{\xi} 
&= \begin{bmatrix} (F + G \theta_1\T) & G \theta_2\T \\ 0 & F \end{bmatrix} \xi +
\begin{bmatrix} 0 \\ G \end{bmatrix} u \\
\label{eq:overp2:xi}
&=: \Fb(\theta) \xi + \Gb u \\
y\0 
&= \theta\T \xi \,.
\end{align}
\end{subequations}
The set of parameters for which the overparameterized model is an input-output equivalent model of the plant is 
\begin{align}
\cP &:= 
\left\{ \theta \in \RR^{2 \nb} ~|~ 
\theta\T ( sI - \Fb(\theta) )\inv \Gb = \frac{\beta(s)}{\al(s)} \right\} \,.
\end{align}
 
\begin{prop}[\cite{ArXivSID26,KREISSEL96}]
\label{prop:thetas}
Suppose Assumptions \ref{as:coprime}-\ref{as:measurements} hold.
Consider \eqref{eq:overp1} with $(F,G)$ defined in \eqref{eq:FG}. Let $\delta(s)$ be any monic polynomial of degree $p = \nb - n$. Select $\theta = (\theta_1,\theta_2) \in \RR^{2 \nb}$ with
\begin{align}
\label{eq:thetas}
\theta_1 = f - \ale \,, \qquad \theta_2 = \be \,,
\end{align}
where $\ale = (\bar{\al}_0, \bar{\al}_1, \ldots, \bar{\al}_{\nb-1})$ and $\be = (\bar{\beta}_0, \bar{\beta}_1, \ldots, \bar{\beta}_{\nb-1})$ are such that
\begin{subequations}
\label{eq:extended}
\begin{align}
\nonumber
\ale(s) 
&:= s^{\nb} + \bar{\al}_{\nb-1} s^{\nb-1} + \cdots + \bar{\al}_1 s + \bar{\al}_0 \\
&= \al(s) \delta(s) \\ 
\be(s) 
&:= \bar{\beta}_{\nb-1} s^{\nb-1} + \cdots + \bar{\beta}_1 s + \bar{\beta}_0 = \beta(s) \delta(s) \,.
\end{align}
\end{subequations}
Then $\theta \in \cP$ if and only if $\theta$ satisfies \eqref{eq:thetas} for some $\delta(s)$.\tqed
\end{prop} 

The previous result captures in precise terms the presence of a pole-zero cancellation in forming the transfer function of the overparameterized model in order for it to match the plant transfer function. One can go further to match the plant and overparameterized model outputs by properly selecting initial conditions of \eqref{eq:overp1} or \eqref{eq:overp2}.

\begin{prop}[\cite{ArXivSID26,KREISSEL96}]
\label{prop:matching}
Suppose Assumptions \ref{as:coprime}-\ref{as:measurements} hold.
Consider \eqref{eq:overp2} with $(F,G)$ defined in \eqref{eq:FG} and $\theta \in \cP$. For each initial condition $\yb(t_0) := (y,\dot{y},\ldots,y^{(n-1)})(t_0)$ of \eqref{eq:theplant}, there exists $\xi(t_0) \in \RR^{2 \nb}$ such that the solution of \eqref{eq:overp2} satisfies: 
$y\0(t) = y(t)$ and $\xi(t) \in \cC(\theta)$, the controllable subspace of \eqref{eq:overp2}, for all $t \ge t_0 \ge 0$. \tqed
\end{prop} 

Let $\xis(t; \theta, \yb(t_0))$ denote the output matching solution in Proposition~\ref{prop:matching} for a fixed $\theta \in \cP$ and initial condition $\yb(t_0) \in \RR^n$.
The next result ensures that we may study the properties of this solution using any $\theta \in \cP$. 

\begin{prop}[\cite{GHADIEH26}]
\label{prop:controllable}
Suppose Assumptions \ref{as:coprime}-\ref{as:measurements} hold.
Consider \eqref{eq:overp2} with $(F,G)$ defined in \eqref{eq:FG}. 
For each initial condition $\yb(t_0) \in \RR^n$ of \eqref{eq:theplant} and for all $\theta, \theta' \in \cP$, we have $\cC(\theta) = \cC(\theta')$ and $\xis(t;\theta,\yb(t_0)) = \xis(t;\theta',\yb(t_0))$ for all $t \ge t_0 \ge 0$.\tqed
\end{prop}

We henceforth drop the dependence on $\theta$ and simply write $\xis(t;\yb(t_0))$ to denote the output matching solution. 
When expedient, we simply write $\xis(t)$.

\subsection{Special Choice of Parameters}

The characterization of $\cP$ in Proposition~\ref{prop:thetas} shows that when $\nb > n$, there exists an infinite number of matching parameters  which are parameterized by a monic polynomial $\delta(s)$ of degree $p = \nb - n$.  If we select $\delta(s) = s^p$, then $\ale = (\Zero_p\,,\al)$ and $\be = (\Zero_p \,, \beta)$ in \eqref{eq:extended}.
Define $\tho := \begin{bmatrix} f\T & \Zero\T_{\nb} \end{bmatrix}\T \in \RR^{2\nb}$. The parameter $\ths \in \cP$ corresponding to this choice of $\delta(s)$ is 
\begin{align} \label{eq:thetastar}
\ths &= \tho +
\begin{bmatrix}
    \Zero_p \\ -\al \\ \Zero_{p} \\ \beta
\end{bmatrix} \,.
\end{align}
It has the useful property of providing a direct readout of the plant parameters $(\al,\beta)$. Our strategy in the sequel will be to develop an adaptive algorithm to recover $\ths$. 

\section{Regressor Redundancy}
In this section we explore the fact that in the overparameterized case $\nb > n$, the non-uniqueness of the parameters in Proposition~\ref{prop:thetas}, namely $\cP$ is not a singleton, manifests itself as a loss of PE in $\xis$. 
This property will ultimately provide a mechanism to identify $\ths$. 
We explore the loss of PE of $\xis$ via a structured form called regressor redundancy. 

\begin{definition}\label{defn:redundant}
    Let $v(t) \in \RR$ and $w(t) \in \RR^q$ be two signals. We say $v$ is \emph{redundant} with respect to (w.r.t.) $w$ if there exists $c \in \RR^q$ such that $v - c\T w$ is not PE. 
    Conversely, $v$ is \emph{not redundant} w.r.t. $w$ if for all $c \in \RR^{q}$, $v - c\T w$ is PE. \tqed
\end{definition}
\begin{example} \label{ex:redundant}
To illustrate these definitions, we consider two cases.
First, let $w(t) = \begin{bmatrix} \sin(t) & \cos(t) \end{bmatrix}\T$ and $v(t) = 2\sin(t) - 3\cos(t) + e^{-t}$. Choosing $c = \begin{bmatrix} 2 & -3 \end{bmatrix}\T$, then $v(t) - c\T w(t) = e^{-t}$, which converges to $0$ as $t \to \infty$. Therefore, $v - c\T w$ is not PE, and $v$ is redundant w.r.t. $w$.
Second, let $w(t) = \sin(t)$ and $v(t) = \sin(2t)$. For any $c \in \RR$, consider $r(t) = \sin(2t) - c\sin(t)$, which contains distinct frequency components. Then for all $c \in \RR$, $v - c\T w$ is PE. Therefore, $v$ is not redundant w.r.t. $w$. 
\tqed
\end{example}

We see from Example~\ref{ex:redundant} that the notion of redundancy captures whether a scalar regressor adds novel excitation to another (vector) regressor. The concept of redundancy also relates to a well-known test for PE \cite[Lemma 6.2]{NARENDRA89}. 

\begin{lemma} \label{lem:stack}
Consider $v(t) \in \RR$, $w(t) \in \RR^q$, and $\wb(t) := (v, w)(t) \in \RR^{q+1}$. Then $\wb$ is PE if and only if $w$ is PE and $v$ is not redundant w.r.t. $w$.
\tqed
\end{lemma}

We now apply our notion of redundancy to $\xis$. 
Let $\xis_{1,i}$ and $\xis_{2,i}$ denote the $i$-th components of $\xis_1$ and $\xis_2$, respectively, for $i = 1, \ldots, \nb$. Define
\begin{align} \label{eq:windowdef}
\xisbi &= \begin{bmatrix}
\xis_{1,i}  & \ldots & \xis_{1,\nb} & \xis_{2,1}  & \ldots & \xis_{2,\nb}
\end{bmatrix}\T \,,
\end{align} 
which contains the components of $\xis_1$ from index $i$ to $\nb$ and the full sub-vector $\xis_2$. 
For each $i = 1, \ldots, \nb$, we test whether $\xis_{1,i}$ is redundant w.r.t. $\xisb$. To that end, we define the {\em residual} 
\begin{align}\label{eq:ri_def}
r^\star_i(t,c_i) := \xis_{1,i}(t) - c_i\T \xisb(t) 
\end{align}
for some $c_i \in \RR^{2\nb-i}$. The following result shows that there is a precise boundary between redundancy and non-redundancy in $\xis$.

\begin{theorem} \label{thm:Redundancy}
Suppose Assumption~\ref{assum:plant} holds. 
Fix any $\yb(t_0) \in \RR^n$ and consider the output matching solution $\xis(t;\yb(t_0))$ of \eqref{eq:overp2}. 
For each $i \in \{ 1, \ldots, \nb \}$ we have 
\begin{enumerate}
    \item[a)] If $i > p$, then $\xis_{1,i}$ is not redundant w.r.t. $\xisb$, and  $\xisbi$ is PE.
    \item[b)] If $i \le p$, then $\xis_{1,i}$ is redundant w.r.t. $\xisb$, $\xisbi$ is not PE, and there exists $c_i^0 \in \RR^{2\nb - i}$ such that $r^\star_i(t, c_i^0)$ converges to zero exponentially. \tqed
\end{enumerate} 
\end{theorem}  
We make three observations. First, if $p = \nb - n = 0$, then case (b) does not occur, so for all $i = 1, \ldots, \nb$, $\xis_{1,i}$ is not redundant w.r.t. $\xisb$, and $\xisbi$ is PE. Second, the full vector $\xis_2$ is always included in the test for redundancy. Since  $(F,G)$ is controllable and $u$ is SR, by~\cite[Lemma 2]{ArXivSID26}, $\xis_2$ is PE. This means that redundancy cannot arise among the components of $\xis_2$. Hence, any loss of PE in $\xis$ arises from redundancy of the components of $\xis_1$.  
Also we note that the redundancy structure revealed in Theorem~\ref{thm:Redundancy} is tied to the specific choice of $(F,G)$ in controllable canonical form. Finally, while Theorem~\ref{thm:Redundancy} guarantees that for each $i \le p$ there exists a \emph{redundancy vector} $c^0_i$ for $\xis$,  this vector may not be unique. The following lemma characterizes the geometric structure of the set of redundancy vectors. 

\begin{lemma} \label{lem:redSubspace}
Suppose the conditions of Theorem~\ref{thm:Redundancy} hold.
Let $1 \le i \le p$ and let $\cWp_{i+1}$ denote the non-PE subspace of $\xisb$.
For any $c \in \RR^{2\nb-i}$, $r^\star_i(t, c)$ is not PE if and only if
\begin{align*}
 c \in c^0_i + \cWp_{i+1} \,.
\end{align*} 
Moreover, $r^\star_i(t,c)$ converges to zero exponentially for every such $c$.
\tqed
\end{lemma}

\section{Adaptive Identification}

We present the key components of our design: an overparameterized observer, a robust redundancy detector, and a parameter adaptation law.

\subsection{Overparameterized Observer}\label{sec:observer}

We consider an {\em overparameterized observer} of the form
\begin{subequations}
\label{eq:observer}
\begin{align}
\dot{\xih}_1 &= F \xih_1 + G y \\
\dot{\xih}_2 &= F \xih_2 + G u \,,
\end{align}
\end{subequations}
with $\xih := (\xih_1, \xih_2) \in \RR^{2 \nb}$; arbitrary initial conditions $\xih(t_0) := (\xih_1(t_0), \xih_2(t_0))$; and $(F,G)$ as in \eqref{eq:FG}.
Let $\xih_{1,i}$ and $\xih_{2,i}$ denote the $i$-th components of 
$\xih_1$ and $\xih_2$ and define $\xih_{i:}$ analogously to $\xis_{i:}$ in \eqref{eq:windowdef}. Also define $\xit := \xih - \xis$. 
Then we have
\begin{align}
\label{eq:xit}
\dot{\xit} 
&= \begin{bmatrix}
F & 0\\
0 & F\\
\end{bmatrix} \xit =: \Fbz \xit  \,.
\end{align}
Thus, $\xih$  converges exponentially to the output matching solution $\xis$ of \eqref{eq:overp2}, even though the latter depends on unknown parameters and initial conditions. 

\subsection{Robust Redundancy Detector}\label{sec:gate}

The second major component of our adaptive design is a redundancy detector that robustly determines whether or not $\xih_{1,i}$ is redundant w.r.t. $\xihb$. To this end, fix $m_p > 0$ and $\veps_p > 0$ and define the convex sets
\begin{align*}
\cS_0 &:= \left\{ \hci \in \RR^{2\nb-i} ~|~  \|\hci\| \le m_p \right\} \,, \\
\cS_1 &:= \left\{ \hci \in \RR^{2\nb-i} ~|~ 
\|\hci\| \le \sqrt{m_p^2 + \veps_p} \right\} \,.
\end{align*}
We use the smooth projection operator of \cite{KRSTIC95}:
\begin{equation}\label{eq:projdef}
\Proj(\eta) :=
\begin{cases}
  \eta \,, & \text{if } \hci \in \interior(\cS_0) \text{ or }\\ 
         &\hci\T \eta \le 0 \,, \\[0.5ex]
  \Big( I - \kappa(\hci) \dfrac{\hci \, \hci\T}{\| \hci \|^{2}} \Big) \eta \,, & \text{if } \hci \in \cS_1 \setminus \interior(\cS_0) \\ & \text{ and } \hci\T \eta > 0 \,,
\end{cases}
\end{equation}
where $\kappa(\hci) := \min\{1,\, \veps_p\inv ( \| \hci \|^2 - m_p^2) \}$. Note that $\kappa = 0$ on $\partial\cS_0$ and $\kappa = 1$ on $\partial\cS_1$, so the outward component of $\eta$ is removed gradually across the boundary layer $\cS_1 \setminus \cS_0$ and fully on $\partial\cS_1$. 

We propose the {\em redundancy detector}
\begin{subequations} 
\label{eq:detector}
\begin{align}
\label{eq:claw}
\dhci &= \Proj\big( \gamma_c\, \rhat\, \xihb \big) \\
\label{eq:rhat}
\rhat & = \xih_{1,i} - \hci\T \xihb  \\
\label{eq:sigr}
\dsigr &= -\veps \sigr + \veps \rhat^2 \,, 
\end{align}   
\end{subequations}
where $i = 1, \ldots, \nb$; $\gamma_c, \veps > 0$ are design parameters; and we select $\sigr(t_0) \ge 0$. 
An interesting feature of this design is that detection of redundancy in \eqref{eq:claw}-\eqref{eq:rhat} has been interpreted as an adaptive control problem. To see this, consider the estimate of the residual in \eqref{eq:rhat}. 
We showed in Section~\ref{sec:observer} that $\xih_{1,i}(t) = \xis_{1,i}(t)$ for $t \ge t_0$, modulo an exponentially vanishing term. 
By Theorem~\ref{thm:Redundancy}(b) and Lemma~\ref{lem:redSubspace}, for $i \le p$ and any $c_i \in c^0_i + \cWp_{i+1}$, we have 
$\xis_{1,i}(t) = c_i\T \xis_{i+1:}(t)$ for $t \ge t_0$, 
again modulo an exponentially vanishing term. Therefore, \eqref{eq:rhat} becomes
\begin{align*}
\rhat &\simeq (c_i - \hci)\T \xihb \,.
\end{align*}
This equation is the static error model of adaptive control \cite[Ch.~7]{NARENDRA89} with $\rhat$ as the error, $c_i$ interpreted as the unknown parameter, and $\xihb$ the known regressor. Due to this interpretation, we utilize a standard gradient law in \eqref{eq:claw} to recover an estimate of $c_i$. The $\Proj$ operator prevents parameter drift along $\cWp_{i+1}$.
The filter \eqref{eq:sigr} is a scalar version of Kreisselmeier's integral algorithm \cite{KREISSEL77} and was used in \cite{UZEDA25,GHADIEH26} to estimate PE subspaces. Here it estimates the level of excitation in $\rhat$. Its relevance is that we cannot use $\rhat$ directly to detect redundancy as $\rhat$ may oscillate and pass through zero, even when it is PE.

To further distinguish redundancy from non-redundancy, we introduce a binary gate function $\rho_i(t)$ with the dwell-time logic of \cite{MORSE96} (see also \cite[Fig.~40]{LIBERZON03}). 
Let $\tau_d > 0$ be a user-specified dwell time and let $t_k := t_0 + k \tau_d$, $k \in \NN$. For each $i \in \{1, \ldots, \nb\}$, the gate is held constant on each interval $[t_k, t_{k+1})$ and is updated at the instants $t_k$, $k \ge 0$, according to
\begin{align}\label{eq:gate}
\rho_i(t) &=
\begin{cases}
0, & \quad \text{if } \sigr(t_k) \leq \stol \,,\\
1, & \quad \text{if } \sigr(t_k) > \stol \,,
\end{cases}
\qquad t \in [t_k, t_{k+1}) \,,
\end{align}
where $\stol > 0$ is a to-be-determined design parameter. Thus, consecutive switches of $\rho_i$ are separated by at least $\tau_d$, which prevents fast switching when $\sigr$ crosses $\stol$ repeatedly.
\subsection{Parameter Adaptation Law}

The final component of the design is an adaptation law to estimate the special parameter \eqref{eq:thetastar}. Define $\psis = \ths - \tho = (\psis_1,\psis_2) \in \RR^{\nb} \times \RR^{\nb}$, its estimate $\psih \in \RR^{2\nb}$, and the parameter estimation error $\psit := \psih - \psis$. Also recall $\xit = \xih - \xis$. The output error is defined as
\begin{align}\label{eq:error}
e &:= (\tho + \psih)\T \xih - y = \psit\T\xih + (\ths)\T\xit \,,
\end{align}
where we have used $y = (\tho + \psis)\T \xis$ by Proposition~\ref{prop:matching}.
Let $I_2$ denote the $2 \times 2$ identity matrix and let $\Lambda := I_2 \otimes \diag(\rho_1, \ldots, \rho_{\nb})$, so that $\rho_i$ acts on the pair $(\psih_{1,i}, \psih_{2,i})$. The proposed parameter adaptation law is
\begin{align}\label{eq:thetalaw}
    \dot{\psih} = -\gamma  \Lambda  \xih  e - \mu (I - \Lambda)  \psih \,,
\end{align}
where $\gamma, \mu > 0$. The first term is the {\em learning term} and the second is the {\em forgetting term} \cite{MEJIA23_SCL,UZEDA23}. By comparing Theorem~\ref{thm:Redundancy} and \eqref{eq:thetastar}, we see that the indices $i \le p$ for which $\xis_{1,i}$ is redundant correspond to  indices when both $\psis_{1,i}$ and $\psis_{2,i}$ are zero. We show in the next section that forgetting is ultimately applied precisely to these components; namely, $(\psih_{1,i},\psih_{2,i})$ with $i \le p$ will be forgotten.

\section{Main Results}
\label{sec:main}

The complete adaptive identification scheme is
\begin{subequations}\label{eq:design}
\begin{align} 
    \label{eq:main_obs1}
    \dot{\xih}_1 & = F \xih_1 + G y \\
    \label{eq:main_obs2}
    \dot{\xih}_2 & = F \xih_2 + G u \\
    \label{eq:main_rhat}
    \rhat        & = \xih_{1,i} - \hci\T \xihb\\
    \label{eq:main_hci}
    \dhci        & = \Proj\big( \gamma_c \rhat \xihb \big) \\
    \label{eq:main_sig}
    \dsigr       & = -\veps \sigr + \veps \rhat^2 \\
    \label{eq:main_rho} 
    \rho_i(t) &=
    \begin{cases}
    0, & \quad \text{if } \sigr(t_k) \leq \stol \,,\\
    1, & \quad \text{if } \sigr(t_k) > \stol \,,
    \end{cases}
    \qquad t \in [t_k, t_{k+1}) \,,\\ 
    \label{eq:main_e}
    e           & = (\tho + \psih)\T \xih - y \\
    \label{eq:main_psi}
    \dot{\psih}  & = -\gamma  \Lambda \xih e - \mu (I - \Lambda) \psih \\
     \label{eq:main_nh}
    \nh(t)       & = \sum_{i=1}^{\nb} \rho_i(t)
\end{align}
\end{subequations}
where $i = 1, \ldots, \nb$; $m_p$, $\stol$, $\veps_p$, $\veps$, $\gamma_c$, $\gamma, \mu > 0$ are design parameters; and $\rho_i$ are updated with dwell time $\tau_d > 0$ as described in Section~\ref{sec:gate}. We have included an estimate of the plant order in the last equation. 
The analysis is conducted in three main steps: (i) we show the observer state $\xih$ preserves the redundancy and PE properties of $\xis$; (ii) we prove that, by proper selection of $\stol$,  the binary gates $\rho_i$ converge to their correct values in finite time; (iii) we prove $\psih$ converges exponentially to $\psis$.

For the first step of the analysis, we show that $\xih$ may be used in place of $\xis$ to detect redundancy. This is because $\xih$ and $\xis$ differ by an exponentially vanishing term, so the PE properties of $\xisbi$ carry over to $\xihbi$ \cite[Lemma~6.5]{NARENDRA89}.
\begin{prop} \label{prop:xihPE}
Suppose Assumption~\ref{assum:plant} holds. Consider the observer \eqref{eq:main_obs1}--\eqref{eq:main_obs2} and let $1 \le i \le \nb$. Then for every $\yb(t_0) \in \RR^n$ and $\xih(t_0) \in \RR^{2\nb}$, $\xihbi$ is bounded, and
\begin{enumerate}
    \item[a)] if $i > p$, then $\xihbi$ is PE;
    \item[b)] if $i \le p$, then $\xihbi$ is not PE. \tqed
\end{enumerate}
\end{prop}
Redundancy is likewise preserved when $\xis$ is replaced by $\xih$.
\begin{prop} \label{prop:xih_redundancy}
Suppose Assumption~\ref{assum:plant} holds. Consider the observer \eqref{eq:main_obs1}--\eqref{eq:main_obs2} and let $1 \le i \le \nb$. Then $\xih_{1,i}$ is redundant w.r.t. $\xihb$ if and only if $\xis_{1,i}$ is redundant w.r.t. $\xisb$. Moreover, for every $c_i \in \RR^{2\nb-i}$, $r_i(t,c_i) = \xih_{1,i} - c\T_i\xihb$ and $r^\star_i(t,c_i) = \xis_{1,i} - c\T_i \xisb$ satisfy 
$r_i(t, c_i) = r^\star_i(t, c_i) + \eps_i(t)$, 
where $\eps_i(t)$ converges to zero exponentially. \tqed
\end{prop} 

\subsection{Correctness of the Redundancy Detector}

The second step of the analysis is to show that \eqref{eq:main_sig}--\eqref{eq:main_rho} classify redundancy correctly.
For the case $i \le p$, we have discussed in Section~\ref{sec:gate} that the redundancy detector aims to identify some element of the set $c^0_i + \cWp_{i+1}$. In order for the gradient algorithm in \eqref{eq:main_hci} to operate correctly with the $\Proj$, we require that $\cS_0$ contain at least one element of this set.
\begin{assumption}\label{assum:proj}
We select $m_p > 0$ such that for all $1 \le i \le p$,
\begin{align*}
\cS_0 \cap  \big( c^0_i + \cWp_{i+1} \big) \neq \emptyset \,,
\end{align*}
where $c^0_i + \cWp_{i+1}$ is given in Lemma~\ref{lem:redSubspace}.
\tqed
\end{assumption} 

The next result presents our main finding on correctness of the redundancy detector. It tells us that the gradient law and its output \eqref{eq:main_hci}-\eqref{eq:main_rhat} are able to correctly classify whether $\rhat$ is PE or not PE, analogous to Theorem~\ref{thm:Redundancy} for $\xis$ and $r^{\star}_i$.
It should be noted that for $i \le p$, this result does not say that $\hci$ converges to a specific value. Instead, the best one can hope for is that $\hci(t)$ asymptotically approaches $c^0_i + \cWp_{i+1}$. 

\begin{theorem} \label{thm:hatrPE}
Suppose Assumptions~\ref{assum:plant} and~\ref{assum:proj} hold. Consider the observer \eqref{eq:main_obs1}--\eqref{eq:main_obs2} and the estimator \eqref{eq:main_rhat}--\eqref{eq:main_hci} with $\hci(t_0) \in \cS_1$ and let $1 \le i \le \nb$. Then for every $\yb(t_0) \in \RR^n$ and $\xih(t_0) \in \RR^{2\nb}$:
\begin{enumerate}
\item[a)] if $i > p$, then $\rhat$ is PE;
\item[b)] if $i \le p$, then $\rhat$ is not PE and $\rhat(t) \to 0$ as $t \to \infty$.
\tqed
\end{enumerate}
\end{theorem}

To make use of Theorem~\ref{thm:hatrPE}, we must be able to distinguish {\em computationally} whether or not $\rhat$ is PE. This is the role of the subspace estimator \eqref{eq:main_sig}, which extracts the level of PE of $\rhat$, information that is passed to the gates that make the final binary decision. Correctness of these two design components depends on some reasonable PE properties of $\rhat$; namely, that if $\rhat$ is PE, its level of PE is not arbitrarily small.
Proposition~\ref{prop:rhatPE} stated next is a purely theoretical result to ensure that for each $\rhat$ with $i > p$, its level of PE can be bounded away from zero. While this level of PE will generally not be known in practice, this quantity is adjustable through the selection of the input magnitude, so a practical algorithm is achievable. In the next result we exploit the fact that the frequencies appearing in $\xihbi$ for $i > p$ are the same frequencies as the input, so a common PE window $T\pe > 0$ can be chosen for all regressors. Thus, let $(\beta_{pe,i}, T\pe)$ be a pair of PE constants for $\xihbi$. 

\begin{prop} \label{prop:rhatPE}
Suppose the conditions of Theorem~\ref{thm:hatrPE} hold and let $i > p$. Then for every $\yb(t_0) \in \RR^n$ and $\xih(t_0) \in \RR^{2\nb}$,
\begin{align}\label{eq:rhatbound}
    \frac{1}{T\pe} \int_t^{t+T\pe} \rhat^2(\tau)\, d\tau \ge \beta^r_i \,, \quad \forall\, t \ge t_0 \,,
\end{align}
where
\begin{align*}
    \beta^r_i := \frac{(\beta_{pe,i})^2}{\|\xihbi\|_{\cLinf}^2 \Big(1 + \sqrt{1 + \gamma_i}\Big)^2} \,, \qquad
    \gamma_i := \frac{1}{2}\, \gamma_c^2\, T\pe^2\, \|\xihbi\|_{\cLinf}^4 \,.
\end{align*}
\tqed
\end{prop} 

We now analyze the subspace estimator \eqref{eq:main_sig}. Its relevant feature is that it admits a steady-state solution $\sigma_i$ that exposes the level of excitation in $\rhat$. Specifically, as proven in \cite[Lemma~2.1]{MARINO22}, $\sigma_i$ is bounded away from zero when $\rhat$ is PE and converges to zero when $\rhat$ does.

\begin{lemma}[\cite{UZEDA25}] \label{lem:sigfilter}
Suppose Assumptions~\ref{assum:plant} and~\ref{assum:proj} hold. Consider the estimator \eqref{eq:main_rhat}-\eqref{eq:main_hci} with $\hci(t_0) \in \cS_1$. Let $\sigma_i$ be the solution of
\begin{align}\label{eq:sigss}
    \dot{\sigma}_i = -\veps \sigma_i + \veps \rhat^2
\end{align}
with initial conditions $\sigma_i(t_0) = 0$ if $i \le p$ and $\sigma_i(t_0) = \veps T\pe \beta^r_i$ if $i > p$, where $T\pe$ and $\beta^r_i$ are as in Proposition~\ref{prop:rhatPE}. 
\begin{enumerate}
    \item[a)] If $i > p$, then $\sigma_i$ is bounded and
    \begin{align*}
        \sigma_i(t) \ge \veps T\pe \beta^r_i e^{-\veps T\pe} \,, \quad \forall\, t \ge t_0 \,.
    \end{align*}
    \item[b)] If $i \le p$, then $\sigma_i(t) \to 0$ as $t \to \infty$.
    \tqed
\end{enumerate}
\end{lemma}
Define $\sigt := \sigr - \sigma_i$ for $i = 1, \ldots, \nb$. Then $\dot{\tilde{\sigma}}_i = -\veps \sigt$ and $\sigt(t) \to 0$ exponentially. This means Lemma~\ref{lem:sigfilter} can be used to  define a threshold to distinguish the redundant ($i \le p$) and not redundant cases ($i > p$), despite $p$ being unknown. Thus, we define
\begin{subequations}
\begin{align} 
\label{eq:sigmamin}
\sigma_{\min} 
&:= \min_{i > p} \veps T\pe \beta^r_i e^{-\veps T\pe} > 0 \\
\stol 
&\in (0, \sigma_{\min}) \,.
\end{align} 
\end{subequations}

Finally, we verify correctness of the gates \eqref{eq:main_rho}.
Let $x(t_0) := (\yb, \xih, \hat{c}_1, \ldots, \hat{c}_{\nb}, \hat{\sigma}_1, \ldots, \hat{\sigma}_{\nb}, \psih)(t_0)$ denote the initial conditions of the plant and of the identification scheme. 
Since we know $\sigr \to \sigma_i$ exponentially, it is  straightforward to show that for $i \le p$, 
there exists a finite time $T_i(x(t_0))$ such that $\sigr(t) \le \stol$ for all $t \ge T_i(x(t_0))$. Similarly, for $i > p$ there exists a finite time $T_i(x(t_0))$ such that $\sigr(t) > \stol$ for all $t \ge T_i(x(t_0))$.  
In both cases, by \eqref{eq:main_rho}, the gate $\rho_i$ is updated within one dwell time of $T_i(x(t_0))$, so it takes its correct value for all $t \ge T_i(x(t_0)) + \tau_d$. Defining
\begin{align*}
    T_\rho(x(t_0)) := \max_{1 \le i \le \nb} T_i(x(t_0)) + \tau_d \,,
\end{align*}
we have shown the following.

\begin{prop}\label{prop:gates}
Suppose Assumptions~\ref{assum:plant} and~\ref{assum:proj} hold. Consider the adaptive identification scheme \eqref{eq:design} with initial conditions $x(t_0)$ satisfying $\hci(t_0) \in \cS_1$ and $\sigr(t_0) \ge 0$. If $0 < \stol < \sigma_{\min}$ and $\tau_d > 0$, then there exists $T_\rho(x(t_0)) \ge t_0$ such that, for all $t \ge T_\rho(x(t_0))$, $\rho_i(t) = 0$ if $i \le p$ and $\rho_i(t) = 1$ if $i > p$. Consequently, $\nh(t) = n$ for all $t \ge T_\rho(x(t_0))$. \tqed 
\end{prop}  

\subsection{Recovery of Plant Parameters}

In the last step of the analysis, we show that the plant order is identified in finite time and the plant parameters are recovered exponentially.

\begin{theorem}\label{thm:paramconv}
Suppose Assumptions~\ref{assum:plant} and~\ref{assum:proj} hold. Consider the adaptive identification scheme \eqref{eq:design} with initial conditions $x(t_0)$ satisfying $\hci(t_0) \in \cS_1$ and $\sigr(t_0) \ge 0$. If $0 < \stol < \sigma_{\min}$ and $\tau_d > 0$, then $\psih(t) \to \psis$ exponentially. \tqed
\end{theorem}

\begin{remark} 
\label{rem:noise}
We make some general remarks on the effect of output noise in our setup. First, $\xih_2$ is computed from the input $u$ generated by the designer, so it is effectively noise-free. In contrast, noise in $y$ may interfere with the separation of redundant and not redundant components of $\xih_1$ via $\stol$. For a given noise level, this is most easily mitigated by increasing the magnitude of $u$, which raises the excitation of the not redundant components of $\xih_1$ without affecting the redundant ones.

A further factor is that for $i \leq p$, $\xihb$ is not PE, so noise may cause $\hci$ to drift along the non-PE subspace $\cWp_{i+1}$ (see Lemma~\ref{lem:redSubspace}), a phenomenon reminiscent of the Rohrs' counterexamples \cite{ROHRS85}. For this reason, we introduced $\Proj$ to bound this drift inside $\cS_1$. Further, $\rhat$ may not vanish in the presence of noise, even in the redundant case.  As a result, $\sigr$ settles at a noise floor that scales with the noise power. Note that robustness of the subspace estimator \eqref{eq:sigr} has already been studied in \cite{UZEDA23,UZEDA25}. Finally, as long as the gates $\rho_i$ converge to the correct values, the exponential stability result in Theorem~\ref{thm:paramconv} guarantees robustness to sufficiently small bounded perturbations \cite[Lemma~9.2]{KHALIL02}.\tqed
\end{remark}

\section{Simulation Example}

We consider an example with $\nb = 4$, $f(s) = (s+1)^{4}$, $\gamma_c = 0.5$, $\veps = 0.01$, $\gamma = 0.5$, $\mu = 0.5$, and $\tau_d = 10$~s. The input is a multisine with $16$ sinusoids whose frequencies are equally spaced on $[0.2, 5]$ rad/s, each of unit amplitude, with random phases. We take $m_p = 6$ and $\veps_p = 3.6$.  It remains to choose $\stol$, which plays the role of the tolerance in a numerical rank test \cite{GOLUB76,STEWART84}. In the noise-free setting, any $\stol \in (0, \sigma_{\min})$ yields the correct classification. Although $\sigma_{\min}$ depends on the initial conditions through $\beta^r_i$, this dependence vanishes asymptotically, since the plant and observer transients decay exponentially. In the presence of noise, $\stol$ also captures the designer's knowledge of the expected noise level distinguishable from relevant information in $\xih$; see \cite{BAKO09}.

In numerical linear algebra, rank tolerances are usually set relative to the scale of the data. Following \cite{GHADIEH26}, we take $\xih_2$ as a reference excitation, since it is driven only by $u$ and is PE by \cite[Lemma~2]{ArXivSID26}. For $j = 1, \ldots, \nb-1$, we run the detector \eqref{eq:main_rhat}--\eqref{eq:main_sig} with $(\xih_{1,i}, \xihb)$ replaced by $(\xih_{2,j}, \xih_{2,j+1:})$, and denote the resulting subspace estimator as $\hat{\sigma}_{u,j}$. Since $\xih_2$ has no redundant components, $\hat{\sigma}_{u,j}$ are eventually bounded away from zero, and we set
\begin{align*}
    \shattol(t) := 0.1 \min_{j} \hat{\sigma}_{u,j}(t) \,.
\end{align*}
Thus, a component of $\xih_1$ is declared redundant when its excitation falls below $10\%$ of the weakest reference excitation level. Finally, the scheme \eqref{eq:design} is initialized at zero and the plant at a random initial condition.

\begin{figure*}[!t]
  \centering
  \begin{subfigure}{0.33\linewidth}
    \includegraphics[width=\linewidth]{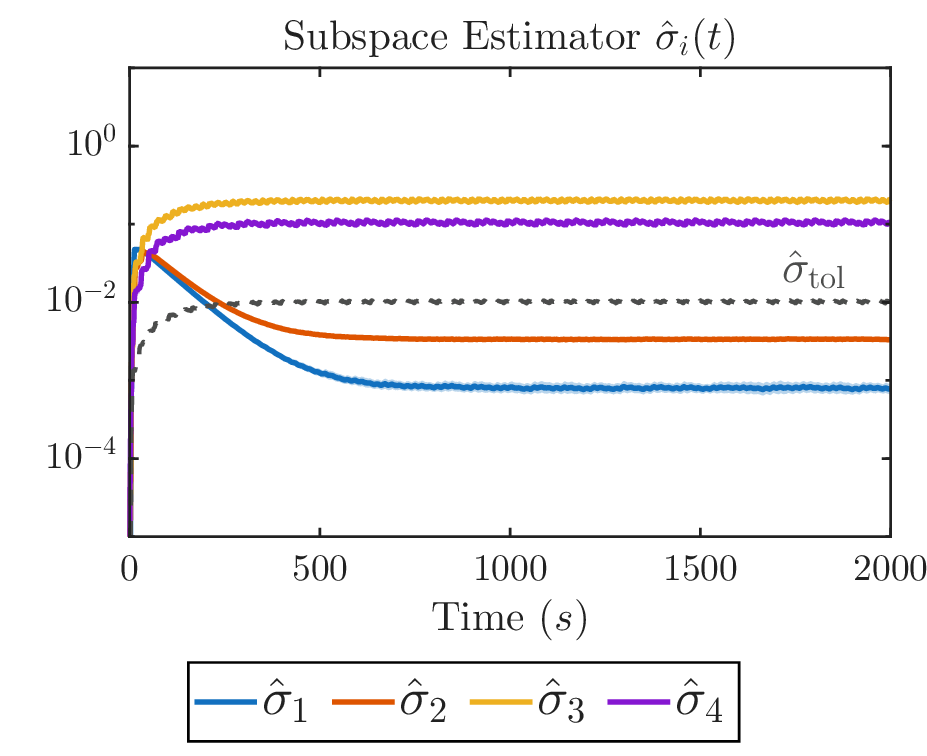}
    \caption{}\label{fig:ex1_sigma}
  \end{subfigure}\hfill
  \begin{subfigure}{0.33\linewidth}
    \includegraphics[width=\linewidth]{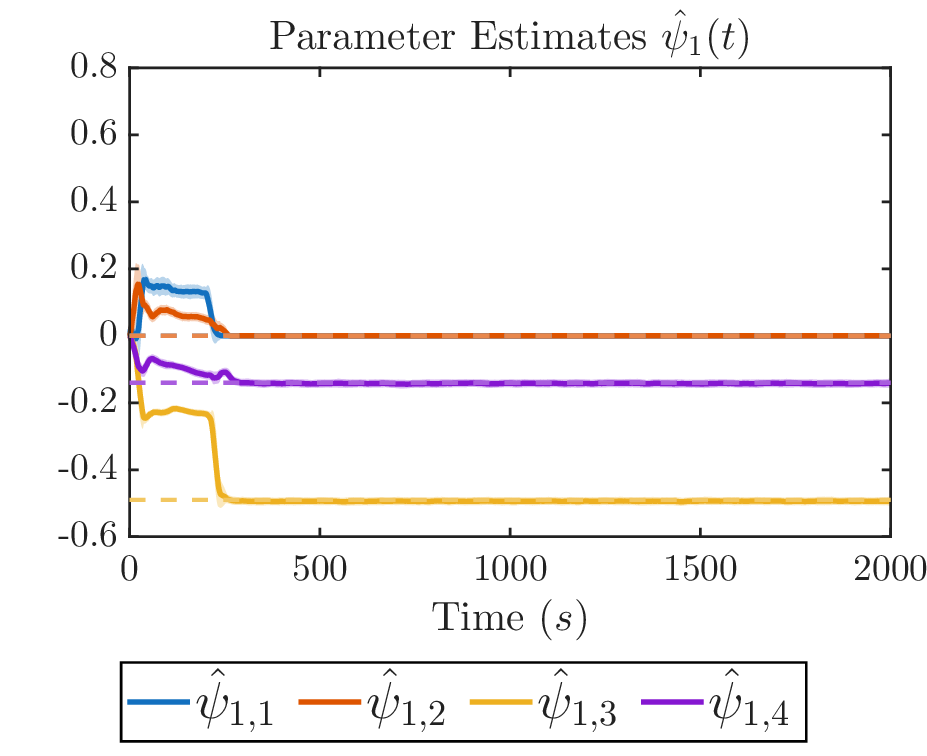}
    \caption{}\label{fig:ex1_psi1}
  \end{subfigure}\hfill
  \begin{subfigure}{0.33\linewidth}
    \includegraphics[width=\linewidth]{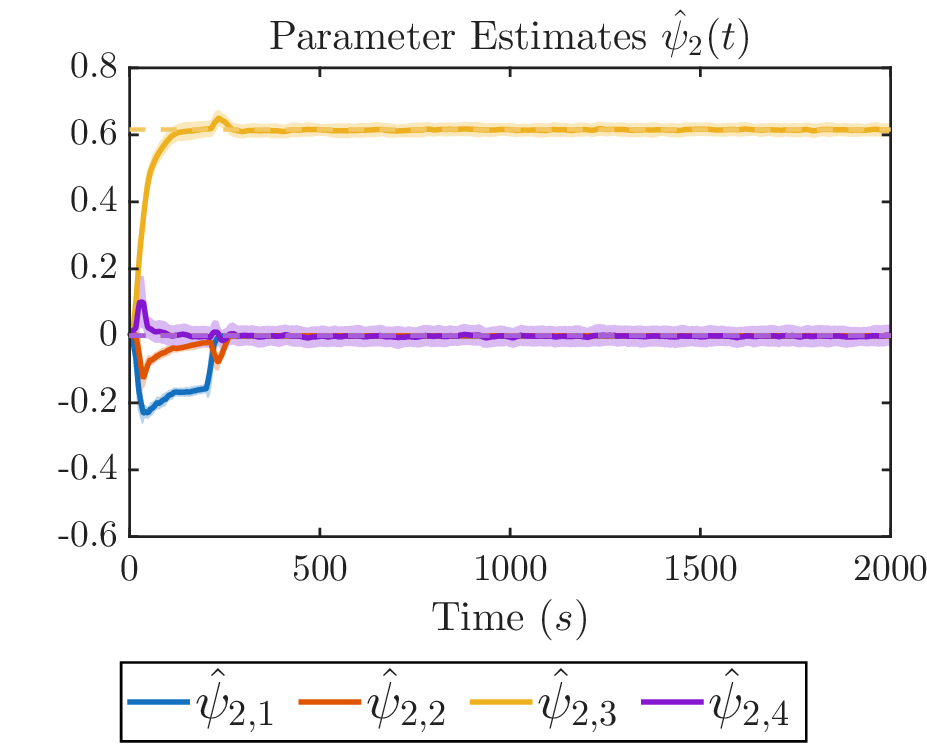}
    \caption{}\label{fig:ex1_psi2}
  \end{subfigure}
  \caption{Adaptive identification experiment with $n = 2$, $\nb = 4$, and a $20$~dB SNR at the output $y$. We plot the mean (solid) and one standard deviation (shaded) over $100$ trials. For $\hat\sigma_i$ the standard deviation is computed on a $\log_{10}\hat\sigma_i$ scale. In (a), we see that $\hat\sigma_1$ and $\hat\sigma_2$ settle below $\shattol$ (dashed), and $\nh(t)$ settles at $n$ in every run. In (b)-(c), $\psih$ converges to a neighborhood of $\psis$ (dashed).}
  \label{fig:ex1}
\end{figure*}

Consider the second-order plant
\begin{align*}
    H(s) = \frac{0.6163}{s^2 + 0.14 s + 0.49}\,,
\end{align*}
for which $\ths$ in \eqref{eq:thetastar} gives
\begin{align*}
    \psis = (0,\, 0,\, -0.49,\, -0.14,\, 0,\, 0,\, 0.6163,\, 0) \,.
\end{align*}
We add zero-mean white Gaussian noise to the measured output $y$ and define the signal-to-noise ratio (SNR) as the ratio of the output signal power to the noise power. We run the scheme \eqref{eq:design} with $t_f = 2000$s for $100$ Monte Carlo trials with different noise realizations. Fig.~\ref{fig:ex1} shows the results at an SNR of $20$~dB. For the redundant components, $\hat\sigma_1$ and $\hat\sigma_2$ settle below $\shattol$, while for the not redundant components,  $\hat\sigma_3$ and $\hat\sigma_4$ settle well above it. The gates converge in every run, and $\nh(t_f) = n$. Once $\rho_1$ and $\rho_2$ switch back to $0$, the entries of $\psih$ switched off by the gates are forgotten and converge to zero, while the remaining entries converge to a neighborhood of the plant coefficients $-\al$ and $\beta$, see Fig.~\ref{fig:ex1_psi1} and Fig.~\ref{fig:ex1_psi2}.

Table~\ref{tab:snr} repeats the experiment at several SNRs. The order is recovered in every run for cases where SNR is $20$~dB or above. In contrast, it is recovered in $30$ trials at $15$~dB and at none for $10$~dB. In the trials where the order is recovered, the coefficient estimates lie in a neighborhood of the true values whose size grows as the SNR decreases.
 
\begin{table}[t]
\centering\footnotesize
\setlength{\tabcolsep}{3pt}
\caption{Model order estimate $\nh$ and parameter estimate $\psih$ at $t_f = 2000$s over 100 trials per SNR. Columns 3-5: mean $\pm$ std of $\psih$ over the trials achieving $\nh(t_f) = n$.}
\label{tab:snr}
\resizebox{\columnwidth}{!}{%
\begin{tabular}{rcccc}
\hline
SNR (dB) & $\nh = n$ & $\psih_{1,3}$ & $\psih_{1,4}$ & $\psih_{2,3}$ \\
\hline
true & & $-0.49$ & $-0.14$ & $0.6163$ \\
\hline
10 & 0/100 & -- & -- & -- \\
15 & 30/100 & $-0.502 \pm 0.039$ & $-0.143 \pm 0.020$ & $0.607 \pm 0.039$ \\
20 & 100/100 & $-0.494 \pm 0.007$ & $-0.142 \pm 0.009$ & $0.616 \pm 0.022$ \\
40 & 100/100 & $-0.490 \pm 0.001$ & $-0.140 \pm 0.002$ & $0.616 \pm 0.002$ \\
$\infty$ & 100/100 & $-0.490$ & $-0.140$ & $0.616$ \\
\hline
\end{tabular}}
\end{table}

\section{Proofs}

\subsection{Proof of Lemma~\ref{lem:kerR}}
Fix $a \in \RR^{q}$ and consider the signal $a\T w(t) \in \RR$. Observe that $a\T w$ is stationary with autocovariance $R[a\T w] = a\T R[w] a$. Because $a\T w$ is a scalar stationary signal, by \cite[Proposition~2.7.1]{SASTRY89}, $a\T w$ is PE if and only if $R[a\T w] > 0$. Conversely, $a\T w$ is not PE if and only if $R[a\T w] = 0$. Since $R[w]$ is symmetric and positive semi-definite, $R[a\T w] = a\T R[w] a = 0$ if and only if $R[w] a = 0$. That is, $a\T w$ is not PE if and only if $a \in \Ker(R[w]) = \cWp$.

\subsection{Proof of Lemma~\ref{lem:stack}} 

First, suppose $\wb$ is PE. Then by \cite[Lemma~6.2]{NARENDRA89}, for all $a \in \RR^{q+1}$, with $a \neq 0$, $a\T \wb$ is PE. In particular, we can choose any $a = (0,a_2) \neq 0$ such that $a\T \wb = a_2\T w$ is PE. Again by \cite[Lemma~6.2]{NARENDRA89}, this implies $w$ is PE. 
Second, choose any $a = (1,-a_2)$. Then $a\T \wb = v - a_2\T w$ is PE, and since $a_2 \in \RR^q$ is arbitrary, $v$ is not redundant w.r.t. $w$.

For the other direction suppose $w$ is PE and $v$ is not redundant w.r.t. $w$. Suppose, by way of contradiction, that $\wb$ is not PE. Then by \cite[Lemma~6.2]{NARENDRA89}, there exists $a := (a_1,a_2) \neq 0$ with $a_1 \in \RR$ and $a_2 \in \RR^q$ such that $a\T \wb$ is not PE. If $a_1 = 0$, then $a_2 \neq 0$ and $a\T \wb = a_2\T w$ is not PE, implying $w$ is not PE by \cite[Lemma~6.2]{NARENDRA89}, a contradiction. Instead, if $a_1 \neq 0$, then $v - a_1\inv a_2\T w =: v - c\T w$ is not PE, which contradicts that $v$ is not redundant w.r.t. $w$. We conclude $\wb$ is PE.

\subsection{Proof of Theorem~\ref{thm:Redundancy}}
Taking the Laplace transform of the output $y(t)$ of the plant \eqref{eq:theplant} with initial condition $\yb(t_0)$ gives
\begin{align}\label{eq:yzero}
  y(s) = \frac{\beta(s)}{\al(s)} u(s) + \frac{\nu(s)}{\al(s)} \,,
\end{align}
where $\nu(s)$ is a polynomial with $\deg (\nu(s)) \le n-1$ whose coefficients are determined by $\yb(t_0)$. Since  $(F,G)$ in \eqref{eq:FG} is in controllable canonical form, 
\begin{align*}
    (sI - F)\inv G = \frac{1}{f(s)} \begin{bmatrix}
    1& s & \ldots & s^{\nb-1} \end{bmatrix}\T \,.
\end{align*}
Therefore, taking the Laplace transform of $\xis$ gives
\begin{equation}\label{eq:xisLaplace}
\begin{split}
    \xis_{1,i}(s) & = \frac{s^{i-1}}{f(s)} y(s) + \frac{\mu_{1,i}(s)}{f(s)} \\
    \xis_{2,i}(s) & = \frac{s^{i-1}}{f(s)} u(s)  + \frac{\mu_{2,i}(s)}{f(s)}
\end{split}
\end{equation}
for $i = 1, \ldots, \nb$. The terms $\mu_{1,i}$ and $\mu_{2,i}$ are polynomials with $\deg(\mu_{1,i}(s)) $, $\deg(\mu_{2,i}(s)) \le \nb - 1$ whose coefficients are determined by $\xis(t_0;\yb(t_0))$. Let $c_{i} = (c_{1,i+1}, \ldots, c_{1,\nb}, c_{2,1}, \ldots, c_{2,\nb}) \in \RR^{2\nb-i}$ be arbitrary and consider the residual $r^\star_i(t,c_i)$ in \eqref{eq:ri_def}. Taking the Laplace transform and using \eqref{eq:yzero} and \eqref{eq:xisLaplace}, we obtain
\begin{align} \label{eq:ri}
\nonumber
r^\star_i(s)
&= \frac{1}{f(s)}
  \Big[ \Big( s^{i-1} - \sum_{j=i+1}^{\nb} c_{1,j} s^{j-1} \Big) y(s)
        \\
        \nonumber
& - \sum_{j=1}^{\nb} c_{2,j} s^{j-1}  u(s)\Big] 
        + \frac{\mu_i(s)}{f(s)} \\
& = \frac{1}{f(s)}
  \Big[ s^{i-1} q(s) \frac{\beta(s)}{\al(s)} + h(s) \Big] u(s) + \eps^\star_i(s) \,,
\end{align}
where
\begin{align} \label{eq:rqh}
     q(s) & := 1 - \sum_{j=i+1}^{\nb} c_{1,j} s^{j-i} \,, \quad
     h(s)   := - \sum_{j=1}^{\nb} c_{2,j} s^{j-1} \,, \\
     \nonumber
     \mu_i(s) & := \mu_{1,i}(s) - \sum_{j=i+1}^{\nb} c_{1,j}\, \mu_{1,j}(s) - \sum_{j=1}^{\nb} c_{2,j}\, \mu_{2,j}(s) \,, \\
     \label{eq:epsLaplace}
     \eps^\star_i(s) & := \frac{s^{i-1} q(s) \nu(s)}{f(s)\al(s)} + \frac{\mu_i(s)}{f(s)} \,.
\end{align}
Notice that $\deg(q(s)) \leq \nb - i$ and $\deg(h(s)) \leq \nb - 1$. The term $\eps^\star_i(s)$  collects the contributions of the initial conditions coming from the plant and $\xis(t_0)$. Both terms in $\eps^\star_i$ are strictly proper with Hurwitz denominators $f(s)\al(s)$ and $f(s)$. Therefore, the Laplace inverse $\eps^\star_i(t)$ converges to zero exponentially. Denote the first term in $r^\star_i(s)$ as $\rnom(s)$, which is of the form 
\begin{align}\label{eq:rnom}
    \nonumber
    \rnom(s) & :=  \frac{1}{f(s)\al(s)} \big[ s^{i-1} q(s)  \beta(s)+ h(s)\al(s) \big] u(s)\\
    &=: T(s) u(s)\,.
\end{align}
Taking the inverse Laplace transform, we obtain
\begin{align}\label{eq:rdecomp}
    r^\star_i(t,c_i) = \rnom(t,c_i) + \eps^\star_i(t) \,.
\end{align}
Since $\eps^\star_i(t) \to 0$ exponentially, \cite[Lemma 6.5]{NARENDRA89} implies that $r^\star_i(t, c_i)$ is PE if and only if $\rnom(t,c_i)$ is PE. Hence the initial conditions $\yb(t_0)$ and $\xis(t_0)$ do not affect the PE properties of $r^\star_i$, and it suffices to analyze $\rnom$.

In what follows we analyze the PE properties of $\rnom$ relative to $i$.  By Assumption~\ref{assum:input}, $u(t)$ is a stationary multisine that is sufficiently rich of order $2\nb$. Suppose $u(t)$ has spectral lines at $\pm \omega_1$, $\ldots$, $\pm\omega_{\nb}$ with $0 < \omega_1 < \cdots < \omega_{\nb}$. The spectral measure $S_u(\omega)$ of $u$ is
\begin{align} \label{eq:Su}
S_u(\omega) &= \sum^{\nb}_{k=1} f_u(\omega_k)\delta(\omega-\omega_k) + \sum^{\nb}_{k=1} f_u(-\omega_k)\delta(\omega+\omega_k) \,,
\end{align}
where $f_u(\omega_k)$, $f_u(-\omega_k) > 0$ and $\delta(\cdot)$ refers to the delta function.
By \cite[Prop~1.6.2]{SASTRY89}, $\rnom$ is stationary and
$S_{\rnom}(\omega)= |T(j\omega)|^2 S_u(\omega)$, where $T(s)$ is defined in \eqref{eq:rnom}. 
By \cite[p.~256]{IOANNOU12}
$R[\rnom] = \frac{1}{2\pi} \int^{+\infty}_{-\infty} S_{\rnom}(\omega) d\omega$. Combined with \eqref{eq:Su}, we obtain 
\begin{align*} 
    R[\rnom] & = \frac{1}{2\pi}\sum^{\nb}_{k= 1} \Big[ |T (j\omega_k)|^2f_u ({\omega_k}) + |T (-j\omega_k)|^2f_u ({-\omega_k})\Big] \,.
\end{align*}
Then $R[\rnom] \ge 0$. Further, since $\rnom$ is scalar and stationary, $\rnom$ is PE if and only if $R[\rnom] > 0$ \cite[Proposition~2.7.1]{SASTRY89}. 
Given the expression above, this is equivalent to the condition that $T(s)$ does not have $2\nb$ distinct zeros at $\pm j\omega_1, \ldots, \pm j\omega_{\nb}$. Let 
\begin{align} \label{eq:Ni}
   N(s) := s^{i-1} q(s) \beta(s) + h(s)\al(s) 
\end{align} 
denote the numerator of $T(s)$. Since the denominator $f(s)\al(s)$ of $T(s)$ is Hurwitz, we find that $\rnom(t)$ is PE if and only if $N(s)$ does not have $2\nb$ distinct zeros at $\pm j\omega_1, \ldots, \pm j\omega_{\nb}$.

Now consider the case $i > p$. Fix $c_i \in \RR^{2\nb-i}$ and suppose, by way of contradiction, that $r^\star_i(t,c_i)$ is not PE. Then $\rnom(t,c_i)$ is not PE, and by the foregoing developments, $N(s)$ has $2\nb$ distinct zeros at $\pm j\omega_1, \ldots, \pm j\omega_{\nb}$. Since $\deg(N(s)) \le 2\nb - 1$, $N(s)$ must be identically zero for all $s \in \CC$. That is,
\begin{align*}
  s^{i-1} q(s) \beta(s) = - h(s) \al(s) \,.
\end{align*}
Since $\al(s)$ and $\beta(s)$ are coprime, $\al(s)$ must divide $s^{i-1}q(s)$. However, $\al(s)$ is Hurwitz, so $\al(s)$ divides $q(s)$. Moreover, $q(0) = 1$ by \eqref{eq:rqh}, so $q(s)$ is not identically zero and $\deg(q(s)) \ge \deg(\al(s)) = n$. However, as noted after \eqref{eq:epsLaplace} and by definition of $p$, $\deg(q(s)) \le \nb - i < \nb - p = n$, a contradiction. Thus, for $i > p$ and for any choice of $c_i \in \RR^{2\nb-i}$, $\rnom(t,c_i)$ and $r^\star_i(t,c_i)$ are PE, and $\xis_{1,i}$ is not redundant w.r.t. $\xisb$.

It remains to show that $\xisbi$ is PE for $i > p$. The pair $(F,G)$ is controllable and $u$ is stationary and SR of order $2\nb \ge \nb$, then $\xis_{\nb+1:} = \xis_2$ is PE by~\cite[Lemma 2]{ArXivSID26}. For each $j = \nb, \ldots, i$ we have $j > p$, so $\xis_{1,j}$ is not redundant w.r.t. $\xis_{j+1:}$. Moreover, if $\xis_{j+1:}$ is PE and $\xis_{1,j}$ is not redundant w.r.t. it, then $\xis_{j:}$ is PE by Lemma~\ref{lem:stack}. Applying this reasoning recursively for $j = \nb, \nb-1, \ldots, i$, we obtain that $\xisbi$ is PE. 

Next, we consider the case where $i \leq p$. Select $c_i \in \RR^{2\nb - i}$ such that $q(s)$ and $h(s)$ in \eqref{eq:rqh} are of the form
\begin{align}\label{eq:qh}
  q(s) = \frac{d_i(s)\al(s)}{\al_0} \,, \qquad
  h(s) = -\, \frac{s^{i-1}d_i(s)\beta(s)}{\al_0} \,,
\end{align}
where $d_i(s)$ is an arbitrary polynomial with $\deg(d_i(s)) \le p - i$, $d_i(0) = 1$, and $\al_0 \neq 0$ because $\al(s)$ is Hurwitz. Substituting these choices of $q(s)$ and $h(s)$ into \eqref{eq:Ni}, we obtain
\begin{align*}
  N(s) = \frac{s^{i-1}d_i(s)\al(s)\beta(s)}{\al_0}
         - \frac{s^{i-1}d_i(s)\beta(s)\al(s)}{\al_0} = 0  
\end{align*}
for all $s \in \CC$. Equivalently, we have $T(s) = 0$ for all $s \in \CC$. 
By \eqref{eq:rnom}, $\rnom(t) = 0$ for all $t \ge t_0$. By \eqref{eq:rdecomp}, $r^{\star}_i(t,c_i)$ vanishes exponentially. Moreover, since $(1,-c_i) \neq 0$
and $(1,-c_i)\T \xisbi = r^\star_i(t,c_i)$ is not PE, then $\xisbi$ is not PE by~\cite[Lemma 6.2]{NARENDRA89}. Choosing $c^0_i := c_i$, we have proved the result.

\subsection{Proof of Lemma~\ref{lem:redSubspace}}

Let $1 \le i \le p$. By Theorem~\ref{thm:Redundancy}(b), there exists $c^0_i \in \RR^{2\nb-i}$ such that $r^\star_i(t,c^0_i) \to 0$ exponentially. Let $c \in \RR^{2\nb-i}$ be arbitrary. Using \eqref{eq:ri_def}, we have
\begin{align} \label{eq:kerxisb}
    r^\star_i(t,c) = r^\star_i(t,c^0_i) + (c^0_i - c)\T \xisb(t) \,.
\end{align}
Since $r^\star_i(t,c^0_i) \to 0$ exponentially, by \cite[Lemma~6.5]{NARENDRA89} it follows $r^\star_i(t,c)$ is PE if and only if $(c^0_i - c)\T \xisb$ is PE. By Assumption~\ref{assum:input}, $\xis$ is the response of a stable LTI system to a stationary multisine plus an exponentially decaying transient, hence $\xis$ and its sub-vectors are stationary. Applying Lemma~\ref{lem:kerR} to $\xisb$, the scalar signal $(c^0_i - c)\T \xisb$ is not PE if and only if $(c^0_i - c) \in \cWp_{i+1} = \Ker(R[\xisb])$. Therefore, $r^\star_i(t,c)$ is not PE if and only if $c \in c^0_i + \cWp_{i+1}$. 

To show exponential convergence of $r^\star_i(t,c)$, let $c \in c^0_i + \cWp_{i+1}$. By the argument above, $r^\star_i(t,c)$ is not PE. Also, $\rnom(t,c)$ is not PE by \eqref{eq:rdecomp} and \cite[Lemma~6.5]{NARENDRA89}. 
Following the proof of Theorem~\ref{thm:Redundancy}, since $\rnom(t,c)$ is not PE, we get $R[\rnom] = 0$, and $N(s)$ has $2\nb$ distinct zeros at $\pm j\omega_1, \ldots, \pm j\omega_{\nb}$. Since $\deg(N(s)) \le 2\nb - 1$, $N(s) = T(s) = 0$ for all $s \in \CC$ and $\rnom(t,c) = 0$ for all $t \ge t_0$. Then by \eqref{eq:rdecomp} we have $r^\star_i(t,c) = \eps^\star_i(t)$, which converges to zero exponentially.

\subsection{Proof of Proposition~\ref{prop:xih_redundancy}}

Let $\xis(t) \equiv \xis(t;\yb(t_0))$ and let $c_i \in \RR^{2\nb-i}$ be arbitrary.  Substituting $\xih(t) = \xis(t) + \xit(t)$ into \eqref{eq:rhat}, we have
\begin{align*}
    r_i(t, c_i) &= \big(\xis_{1,i}(t) + \xit_{1,i}(t)\big) - c_i\T \big(\xisb(t) + \xibt(t)\big) \\
    &= \big(\xis_{1,i}(t) - c_i\T \xisb(t)\big) + \big(\xit_{1,i}(t) - c_i\T \xibt(t)\big) \\
    &= r^\star_i(t, c_i) + \eps_i(t) \,,
\end{align*}
where $r^\star_i(t, c_i)$ is given in \eqref{eq:ri_def}, and $\eps_i(t) := \xit_{1,i}(t) - c_i\T \xibt(t)$. Since $\xit(t) \to 0$ exponentially by \eqref{eq:xit}, also $\eps_i(t) \to 0$ exponentially. By~\cite[Lemma 6.5]{NARENDRA89}, the addition of a decaying term does not alter the PE property of a signal. Therefore, $r_i(t, c_i)$ is PE if and only if $r^\star_i(t, c_i)$ is PE. It follows from Definition~\ref{defn:redundant} that the redundancy properties of  $\xih$ identically follow those of $\xis$.

\subsection{Proof of Theorem~\ref{thm:hatrPE}}

Consider the case where $i > p$. First, because $\Proj(\eta)$ satisfies \eqref{eq:projdef}, we may write
\begin{align*}
\Proj(\gamma_c \rhat \xihb) = \gamma_c \rhat P_i(\hci) \xihb
\end{align*}
where
\begin{align}\label{eq:qi}
    P_i(\hci) :=
    \begin{cases}
        I \,, & \text{if } \hci \in \interior(\cS_0) \text{ or} \\
              & \rhat\, \hci\T \xihb \le 0 \,, \\[0.5ex]
        I - \kappa(\hci)\dfrac{\hci \hci\T}{\|\hci\|^2} \,, & \text{otherwise} \,.
    \end{cases}
\end{align}
In the first case of \eqref{eq:qi}, $\| P_i(\hci) \| = 1$. In the second case, $\|\hci\| \ge m_p > 0$ and $\kappa = \kappa(\hci) \in [0,1]$, so for any $v \in \RR^{2\nb -i}$,
\begin{align*}
    \| P_i(\hci) v \|^2
    &= \|v\|^2 - 2\kappa\,\frac{(\hci\T v)^2}{\|\hci\|^2} + \kappa^2\frac{(\hci\T v)^2}{\|\hci\|^2} \\
    &= \|v\|^2 - \kappa(2-\kappa)\,\frac{(\hci\T v)^2}{\|\hci\|^2} \le \|v\|^2 \,.
\end{align*}
Therefore, $\| P_i(\hci) \| \le 1$.

Next, we introduce the notation $\hat{C}_i := (1,-\hat{c}_i) \in \RR^{2 \nb - i + 1}$. Notice by construction $\| \hat{C}_i \| \ge 1$.
Now \eqref{eq:detector} can be rewritten as
\begin{subequations} 
\label{eq:detector2}
\begin{align}
\dot{\hat{C}}_i &= -\gamma_c \begin{bmatrix}
        0 \\ P_i(\hci) \xihb
    \end{bmatrix} \rhat =: K_i(t,\hci) \xihbi\T(t) \hat{C}_i \\
\rhat & = \xihbi\T(t) \hat{C}_i \,.
\end{align}   
\end{subequations}
A standard approach is to transform \eqref{eq:detector2} to a linear time-varying (LTV) system. That is, we evaluate $K_i(t, \hci)$ along a specific solution $\hci(t)$ such that there is one LTV system associated with each initial condition of \eqref{eq:detector2} \cite{LORIA02,BERNARD19}. Then we study the PE properties of the solution of the LTV system, which then carry over to the equivalent solution of \eqref{eq:detector2}.

To this end, fix $\yb(t_0) \in \RR^n$, $\xih(t_0) \in \RR^{2\nb}$, $\hat{C}_i(t_0) = (1,-\hci(t_0))$, and $\hci(t_0) \in \cS_1$. Let $\hci(t)$ be the corresponding solution of \eqref{eq:claw}. Along this solution, $K_i(t) := K_i(t,\hci(t))$ is a fixed function of time. Now we have an LTV system 
\begin{subequations} \label{eq:nu}
\begin{align}
\dot{\nu} &= K_i(t)  \xihbi\T(t) \nu \\ 
y_\nu     &= \xihbi\T(t) \nu \,.
\end{align}
\end{subequations}
To analyze \eqref{eq:nu}, we first consider \eqref{eq:nu} without output injection \begin{subequations}\label{eq:nu_zero} \begin{align}
\dot{\nu} &= 0 \\ 
y_\nu     &= \xihbi\T(t) \nu \,.
\end{align}
\end{subequations}
By Proposition~\ref{prop:xihPE}, $\xihbi$ is bounded and PE. Let $(\beta_{pe,i}, T\pe)$ denote its PE constants. Then it follows that
\begin{align*}
    \beta_{pe,i} T\pe I \preceq \int_t^{t+T\pe} \xihbi(\tau) \xihbi\T(\tau)\, d\tau \preceq  \|\xihbi\|_{\cLinf}^2 T\pe I \,, \forall t \ge t_0 \,.
\end{align*}
Thus, the pair $(\xihbi\T \,, 0)$ is uniformly completely observable (UCO). The pair $(\xihbi\T, K_i\xihbi\T)$ of \eqref{eq:nu} is obtained from $(\xihbi\T, 0)$ by the output injection $K_i(t)\, y_\nu$. Notice that 
\begin{align}\label{eq:Kbound}
\|K_i(t,\hci)\| 
\nonumber
&\le \gamma_c \| P_i(\hci) \| \| \xihb(t) \| \\
&\le \gamma_c \| \xihb \|_{\cLinf} \le \gamma_c \| \xihbi \|_{\cLinf} \,, 
\end{align}
since $\xihb$ is a sub-vector of $\xihbi$. By \cite[Lemma~2]{ZHANG15}, UCO is preserved under  bounded output injection, so $(\xihbi\T, K_i\xihbi\T)$ is UCO and there exists $\beta^r_i
(\beta_{pe,i}, T\pe, \|\xihbi\|_{\cLinf}, \gamma_c) 
> 0$ such that every solution $\nu$ of \eqref{eq:nu} satisfies
\begin{align}\label{eq:rhatPE}
      \int_t^{t+T\pe} |\xihbi\T(\tau) \nu(\tau)|^2\, d\tau \geq T\pe \beta^r_i \|\nu(t)\|^2 \,, \quad \forall\, t \ge t_0 \,.
\end{align}
When the initial conditions of the LTV \eqref{eq:nu} and \eqref{eq:detector2} are equal, $\hat{C}_i(t) = \nu(t)$ for all $t \ge t_0$. Therefore, we can apply \eqref{eq:rhatPE} to $\hat{C}_i$, and using $\|\hat{C}_i(t)\| \ge 1$ we obtain
\begin{align*}
   \int_t^{t+T\pe} \rhat^2(\tau)\, d\tau \ge T\pe \beta^r_i \|\hat{C}_i(t)\|^2  \ge T\pe \beta^r_i \,, \quad \forall t \ge t_0 \,.
\end{align*}
Moreover, $\hci(t) \in \cS_1$ for all $t \ge t_0$ by \cite[Appendix~E]{KRSTIC95} and $\xihbi$ is bounded, so $\rhat$ is bounded and PE. 
 
Next, consider the case when $i \le p$. By Assumption~\ref{assum:proj}, there exists $c_i \in \cS_0 \cap \big( c^0_i + \cWp_{i+1} \big)$, and by Lemma~\ref{lem:redSubspace}, $r^\star_i(t, c_i)$ converges to zero exponentially. By Proposition~\ref{prop:xih_redundancy},
\begin{align*}
    r_i(t, c_i) = \xih_{1,i}(t) - c_i\T \xihb(t) = r^\star_i(t, c_i) + \eps_i(t) =: \eps(t) \,,
\end{align*}
where $\eps_i(t) = \xit_{1,i}(t) - c_i\T \xibt(t)$ captures the observer transients. Since $r^\star_i(t,c_i)$ and $\eps_i(t)$ both converge to zero exponentially, then $\eps(t) \to 0$ exponentially.

Define the estimation error $\ct = \hci - c_i$. Then $\rhat(t)$ can be rewritten as:
\begin{align}
\nonumber
 \rhat(t) & = \xih_{1,i}(t) - \hci\T(t) \xihb(t)\\
 \label{eq:rhieps}
          & = \eps(t) - \ct\T(t) \xihb(t) \,.
\end{align}
Consider the Lyapunov function candidate $V(\ct) = \frac{1}{2\gamma_c} \ct\T \ct$. Taking the time derivative of $V$ along solutions of \eqref{eq:claw}, we obtain:
\begin{align*}
\dot{V} &= \frac{1}{\gamma_c} \ct\T \dhci = \frac{1}{\gamma_c} \ct\T \Proj\big( \gamma_c \rhat \xihb \big) 
\end{align*}
Since $\hci(t) \in \cS_1$ for all $t \ge t_0$ by~\cite[Appendix~E]{KRSTIC95} and $c_i \in \cS_0$, we obtain
\begin{align*}
    (\hci - c_i)\T \Proj\big( \gamma_c \rhat \xihb \big) \le (\hci - c_i)\T \big( \gamma_c \rhat \xihb \big)\,.
\end{align*} 
This yields
\begin{align*}
    \dot{V} \le \ct\T \rhat \xihb = \rhat \big( \ct\T \xihb \big)\,.
\end{align*}
Substituting $\ct\T \xihb = \eps - \rhat$ from \eqref{eq:rhieps}, we have
\begin{align*}
\dot{V} & \le \rhat(\eps - \rhat) = -\rhat^2 + \rhat\eps\,.
\end{align*}
Applying Young's inequality, $\rhat\eps \le \frac{1}{2}\rhat^2 + \frac{1}{2}\eps^2$, the bound becomes
\begin{align*}
\dot{V} \le -\frac{1}{2}\rhat^2 + \frac{1}{2}\eps^2 \,.
\end{align*}
Integrating both sides from $t_0$ to $t$ gives
\begin{align*}
V(t) - V(t_0) \le -\frac{1}{2} \int_{t_0}^t \rhat(\tau)^2 d\tau + \frac{1}{2} \int_{t_0}^t \eps(\tau)^2 d\tau \,.
\end{align*}
Since $V(t) \ge 0$, rearranging the terms yields:
\begin{align*}
\int_{t_0}^t \rhat(\tau)^2 d\tau \le 2V(t_0) + \int_{t_0}^t \eps(\tau)^2 d\tau \,.
\end{align*}
Since $\eps(t)$ converges to zero exponentially, $\eps \in \cL_2$ and the integral $\int_{t_0}^\infty \eps(\tau)^2 d\tau$ is finite. Therefore, $\int_{t_0}^\infty\rhat(\tau)^2 d\tau < \infty$, which implies $\rhat \in \cL_2$. 

Next, we establish that $\rhat$ and $\dot{\hat{r}}_i$ are bounded. Since $u$ is bounded and $\al(s)$ is Hurwitz, $y$ is bounded. Since $F$ is Hurwitz, \eqref{eq:observer} implies that $\xih$ and $\dot{\xih}$ are bounded. By~\cite[Appendix~E]{KRSTIC95}, $\hci(t) \in \cS_1$ for all $t \ge t_0$, so $\hci$ is bounded, and therefore $\rhat = \xih_{1,i} - \hci\T \xihb$ is bounded. By \cite[Appendix~E]{KRSTIC95}, $\|\dhci\| \le \gamma_c |\rhat|\, \|\xihb\|$, so $\dhci$ is bounded. Consequently,
\begin{align*}
    \dot{\hat{r}}_i = \dot{\xih}_{1,i} - \dhci\T \xihb - \hci\T \dot{\xih}_{i+1:}
\end{align*}
is bounded, that is, $\dot{\hat{r}}_i \in \cLinf$. Thus, $\rhat(t)$ is uniformly continuous. Since $\rhat \in \cL_2 \cap \cLinf$ and $\rhat$ is uniformly continuous, the $\cL_2$ form of Barbalat's lemma in \cite{TAO97,SUN23} gives $\lim_{t \to \infty} \rhat(t) = 0$. Thus, $\rhat(t)$ is not PE. 

\subsection{Proof of Proposition~\ref{prop:rhatPE}} 
We follow the proof of \cite[Lemma~2]{ZHANG15}, applied to the linear time-varying system \eqref{eq:nu}, with $A = 0$, $C = \xihbi\T$, and $K = -K_i$. Let $i > p$ and let $\tilde{\Phi}_i(\tau, t)$ be the state transition matrix of \eqref{eq:nu}. Fix $t \ge t_0$. Throughout, let $\tau \in [t, t + T\pe]$. The state transition matrix of \eqref{eq:nu} satisfies 
\begin{align}\label{eq:zz17}
    \tilde{\Phi}_i(\tau, t) - I = \int_t^\tau K_i(\tau') \xihbi\T(\tau') \tilde{\Phi}_i(\tau', t) d\tau' .
\end{align}
Let $v \in \RR^{2\nb - i + 1}$ be an arbitrary unit vector. Multiplying \eqref{eq:zz17} on the left by $\xihbi\T(\tau)$ and on the right by $v$, and integrating the square of both sides over the window,
\begin{align*}
    \int_t^{t+T\pe} & \Big(\xihbi\T(\tau) \tilde{\Phi}_i(\tau, t) v - \xihbi\T(\tau) v \Big)^2 d\tau\\
    & = \int_t^{t+T\pe} \Big( \xihbi\T(\tau) \int_t^\tau K_i(\tau') \xihbi\T(\tau') \tilde{\Phi}_i(\tau', t) v \, d\tau' \Big)^2 d\tau .
\end{align*}
Expanding the square on the left-hand side and rearranging,
\begin{align}
    &\int_t^{t+T\pe} \big(\xihbi\T(\tau) v\big)^2 d\tau
     = 2\int_t^{t+T\pe} \big(\xihbi\T(\tau) v\big) \big(\xihbi\T(\tau) \tilde{\Phi}_i(\tau, t) v\big) d\tau \nonumber \\
    & -\chi_i + \int_t^{t+T\pe} \Big( \xihbi\T(\tau) \int_t^\tau K_i(\tau') \xihbi\T(\tau') \tilde{\Phi}_i(\tau', t) v \, d\tau' \Big)^2 d\tau \,,\label{eq:zz18}
\end{align}
where
\begin{align}\label{eq:zz20}
    \chi_i & := \int_t^{t+T\pe} \big(\xihbi\T(\tau) \tilde{\Phi}_i(\tau, t) v\big)^2 d\tau = v\T M_i(t) v \\
    \nonumber
    M_i(t) & := \int_t^{t+T\pe} \tilde{\Phi}_i\T(\tau, t) \xihbi(\tau) \xihbi\T(\tau) \tilde{\Phi}_i(\tau, t) d\tau .
\end{align}
We bound each term of \eqref{eq:zz18}. By Proposition~\ref{prop:xihPE}, $\xihbi$ is PE with constants $(\beta_{pe,i}, T\pe)$, and it is bounded. Since $\|v\| = 1$, 
\begin{align}\label{eq:zz19}
    T\pe \beta_{pe,i} \le \int_t^{t+T\pe} \big(\xihbi\T(\tau) v\big)^2 d\tau \le T\pe \|\xihbi\|_{\cLinf}^2 .
\end{align}
By the Cauchy--Schwarz inequality, 
\begin{align*}
    &2\int_t^{t+T\pe} \big(\xihbi\T(\tau) v\big) \big(\xihbi\T(\tau) \tilde{\Phi}_i(\tau, t) v\big) d\tau \\
    & \le 2 \left|\int_t^{t+T\pe} \big(\xihbi\T(\tau) v\big) \big(\xihbi\T(\tau) \tilde{\Phi}_i(\tau, t) v\big) d\tau\right|\\
    &\le 2 \sqrt{\int_t^{t+T\pe} (\xihbi\T(\tau) v)^2 d\tau}\sqrt{\int_t^{t+T\pe}(\xihbi\T(\tau)\tilde{\Phi}_i(\tau, t) v)^2 d\tau} \,.
\end{align*}
Therefore, using \eqref{eq:zz19} and \eqref{eq:zz20}, we obtain
\begin{align*}
    2\int_t^{t+T\pe} \big(\xihbi\T(\tau) v\big) \big(\xihbi\T(\tau) \tilde{\Phi}_i(\tau, t) v\big) d\tau \le 2  \|\xihbi\|_{\cLinf} \sqrt{T\pe \chi_i} .
\end{align*}
Finally, since $\|K_i(t)\| \le \gamma_c \|\xihbi\|_{\cLinf}$ and using the Cauchy--Schwarz inequality and \eqref{eq:zz20}, we obtain
\begin{align*}
    \Big( \xihbi\T&(\tau) \int_t^\tau K_i(\tau') \xihbi\T(\tau') \tilde{\Phi}_i(\tau', t) v d\tau' \Big)^2\\
    & \le \|\xihbi\|_{\cLinf}^2 \int_t^\tau \|K_i(\tau')\|^2 d\tau' \int_t^\tau \big(\xihbi\T(\tau') \tilde{\Phi}_i(\tau', t) v\big)^2 d\tau' \\
    & \le \gamma_c^2 \|\xihbi\|_{\cLinf}^4 (\tau - t) \chi_i \,.
\end{align*}
Integrating over $\tau \in [t, t+T\pe]$ yields
\begin{align*}
    \int_t^{t+T\pe} \Big( \xihbi\T(\tau) \int_t^\tau K_i(\tau')  \xihbi\T(\tau') \tilde{\Phi}_i(\tau', t) v \, d\tau' \Big)^2 d\tau \le \gamma_i \chi_i
\end{align*}
where $\gamma_i := \frac{1}{2} \gamma_c^2 T\pe^2 \|\xihbi\|_{\cLinf}^4$. Therefore, using \eqref{eq:zz18} and \eqref{eq:zz19}, and grouping the bounds above, we have
\begin{align}\label{eq:zz26}
    T\pe \beta_{pe,i} \le -\chi_i + 2 \|\xihbi\|_{\cLinf} \sqrt{T\pe \chi_i} + \gamma_i \chi_i \,.
\end{align}
As in \cite{ZHANG15}, we treat \eqref{eq:zz26} as a polynomial inequality. By \eqref{eq:zz20}, $\chi_i \geq 0$. Moreover, $\chi_i \neq 0$ since \eqref{eq:zz26} does not hold at $\chi_i = 0$. Thus, let $\chi_i' := \sqrt{\chi_i} > 0$.  Then \eqref{eq:zz26} becomes
\begin{align*}
    0 \le (\gamma_i - 1)(\chi'_i)^2 + 2 \|\xihbi\|_{\cLinf} \sqrt{T\pe}\, \chi'_i - T\pe \beta_{pe,i} \,.
\end{align*} 
If $\gamma_i \ne 1$,  the right-hand side is a quadratic in $\chi'_i$ with roots at
\begin{align*}
    \frac{\sqrt{T\pe}\,\Big(-\|\xihbi\|_{\cLinf} \pm \sqrt{\|\xihbi\|_{\cLinf}^2 + (\gamma_i - 1)\beta_{pe,i}}\Big)}{\gamma_i - 1} \,.
\end{align*}
Notice that for both $\gamma_i > 1$ and $\gamma_i < 1$, the root with the plus sign is the smallest positive root, and $\chi'_i$ is at least this root. Multiplying it by its conjugate, we obtain
\begin{align}\label{eq:chibound}
    \chi'_i \ge \frac{\sqrt{T\pe}\,\beta_{pe,i}}{\|\xihbi\|_{\cLinf} + \sqrt{\|\xihbi\|_{\cLinf}^2 + (\gamma_i - 1)\beta_{pe,i}}} \,.
\end{align}
If $\gamma_i = 1$, the inequality is linear in $\chi'_i$ and gives 
\begin{align*}
    \chi'_i \ge \frac{\sqrt{T\pe} \beta_{pe,i} }{2\|\xihbi\|_{\cLinf}}
\end{align*}
which is \eqref{eq:chibound} with $\gamma_i = 1$. Hence \eqref{eq:chibound} holds in all cases.
By \eqref{eq:zz19}, $\beta_{pe,i} \le \|\xihbi\|_{\cLinf}^2$, so $(\gamma_i - 1)\beta_{pe,i} \le \gamma_i \|\xihbi\|_{\cLinf}^2$ and
\begin{align*}
    \sqrt{\|\xihbi\|_{\cLinf}^2 + (\gamma_i - 1)\beta_{pe,i}} \le \|\xihbi\|_{\cLinf}\sqrt{1 + \gamma_i} \,.
\end{align*}
Squaring the bound on $\chi'_i$ and using this inequality,
\begin{align*}
    v\T M_i(t) v = (\chi'_i)^2 \ge \frac{T\pe (\beta_{pe,i})^2}{\|\xihbi\|_{\cLinf}^2 \Big(1 + \sqrt{1 + \gamma_i}\Big)^2} =: T\pe \beta^r_i \,.
\end{align*}
Since the unit vector $v$ is arbitrary, $M_i(t) \succeq T\pe \beta^r_i I$.

Finally, $\hat{C}_i$ is the solution of \eqref{eq:nu} with $\nu(t) = \hat{C}_i(t)$, so $\hat{C}_i(\tau) = \tilde{\Phi}_i(\tau, t) \hat{C}_i(t)$ and, by \eqref{eq:detector2}, $\rhat(\tau) = \xihbi\T(\tau) \tilde{\Phi}_i(\tau, t) \hat{C}_i(t)$. Hence, using $\|\hat{C}_i(t)\| \ge 1$,
\begin{align*}
    \frac{1}{T\pe} \int_t^{t+T\pe} \rhat^2(\tau) d\tau = \frac{1}{T\pe} \hat{C}_i\T(t) M_i(t) \hat{C}_i(t) & \ge \beta^r_i \|\hat{C}_i(t)\|^2\\
    &\ge \beta^r_i 
\end{align*} 
which completes the proof for all $t\ge t_0$. 

\subsection{Proof of Theorem~\ref{thm:paramconv}} 
Fix the initial condition $x(t_0)$ and let $T_\rho = T_\rho(x(t_0))$ be the finite time given by Proposition~\ref{prop:gates}. From \eqref{eq:xit}, \eqref{eq:error}, and \eqref{eq:thetalaw}, the error dynamics are
\begin{subequations}\label{eq:errsys}
\begin{align}\label{eq:errsys1}
    \dot{\psit} &= -\gamma\Lambda\xih\big(\xih\T\psit + (\ths)\T\xit\big) - \mu(I - \Lambda)(\psit + \psis) \\
    \dot{\xit} &= \Fbz\xit \,.
\end{align}
\end{subequations} 

The proof is composed of three steps. First, we show that $\psit$ remains bounded on $[t_0, T_\rho]$, where the gates may be incorrect. Second, we show that $(\psit, \xit)$ converges to zero exponentially on $[T_\rho, \infty)$, where the gates are correct. Third, we combine the two bounds.

\emph{Step 1: $t \in [t_0, T_\rho]$.} \
On this interval $\Lambda(t)$ need not correctly separate redundant and not redundant components of $\xih_1$. Nevertheless, by the dwell-time requirement on \eqref{eq:gate}, each $\rho_i$ switches finitely many times on the bounded interval $[t_0, T_\rho]$, so $\Lambda(t)$ is piecewise constant with finitely many discontinuities. Hence, by \cite[Theorem~3.2]{KHALIL02}, the solution $\psit$ exists and is continuous on $[t_0, T_\rho]$. This interval is compact since $T_\rho$ is finite by Proposition~\ref{prop:gates}. Hence there exists $k_1 > 0$ such that
\begin{align}\label{eq:firstint}
    \|\psit(t)\| \le k_1 \,, \quad \forall\, t \in [t_0, T_\rho] \,.
\end{align}
\emph{Step 2: $t \ge T_\rho$.} \ 
By Proposition~\ref{prop:gates}, the gates $\rho_i$, $i = 1, \ldots, \nb$, are constant and the parameter adaptation law splits into a pure forgetting law and a pure gradient law:
\begin{align} \label{eq:split}
   \begin{bmatrix}
        \dot{\psih}_{1,i}\\
         \dot{\psih}_{2,i}
   \end{bmatrix} &= \begin{cases}
       -\mu\begin{bmatrix}
        \psih_{1,i}\\
        \psih_{2,i}
   \end{bmatrix} \,, & \qquad i \le p \,, \\
   -\gamma \begin{bmatrix}
        \xih_{1,i}\\
        \xih_{2,i}
    \end{bmatrix} e \,, & \qquad i > p \,.
   \end{cases} 
\end{align}
Let $\psit_R \in \RR^{2p}$ collect the pairs $(\psit_{1,i}, \psit_{2,i})$ at the redundant levels $i \le p$, and $\psit_N \in \RR^{2n}$ the pairs at the not redundant levels $i > p$. 
Using the same indices, we split $\xih$, $\xis$, and $\xit$ to generate 
$(\xih_R,\xih_N)$, $(\xis_R,\xis_N)$, and $(\xit_R,\xit_N)$ such that
\begin{align*}
\xih\T\psit &= \xih_N\T\psit_N + \xih_R\T\psit_R \\
e           &= \xih_N\T\psit_N + \xih_R\T\psit_R + (\tho + \psis)\T\xit \,.
\end{align*}
Note that by construction due to \eqref{eq:thetastar}, $\psit_R = \psih_R$. Finally, we obtain the expanded error dynamics 
\begin{subequations}
\label{eq:psitilde}
\begin{align}
\label{eq:psiRtilde}
\dot{\psit}_R &= -\mu \psit_R \\
\label{eq:psiNtilde}
\dot{\psit}_N 
&= -\gamma \xis_N (\xis_N)\T \psit_N + g(t, \psit_N, \psit_R, \xit) \\
\dot{\xit} &= \Fbz\xit \,,
\end{align}
\end{subequations}
where
\begin{align*}
g(t, &\psit_N, \psit_R, \xit) := -\gamma\big(\xit_N(\xis_N)\T 
+ \xis_N\xit_N\T + \xit_N\xit_N\T\big)\psit_N \\
&- \gamma \big(\xis_N + \xit_N\big)\Big(\big(\xis_R + \xit_R\big)\T\psit_R 
 + (\tho + \psis)\T\xit\Big) \,.
\end{align*}

The vector $\xih_N$ is a sub-vector of $\xih_{p+1:}$, which is PE by Proposition~\ref{prop:xihPE}. Hence $\xih_N$ is PE. Since $\xis_N = \xih_N - \xit_N$ and $\xit_N \to 0$ exponentially, $\xis_N$ is PE by \cite[Lemma~6.5]{NARENDRA89}. By \cite[Theorem~2.5.1]{SASTRY89}, the equilibrium $\psit_N = 0$ of the nominal system
\begin{align}\label{eq:nominal}
    \dot{\psit}_N = -\gamma \xis_N (\xis_N)\T \psit_N
\end{align}
is globally exponentially stable. By the global version of \cite[Theorem~4.14]{KHALIL02} there exist $V_0 : [t_0,\infty) \times \RR^{2n} \to \RR$ and $a_1, a_2, a_3, a_4 > 0$ such that
\begin{subequations}\label{eq:converse}
\begin{align}
    a_1\|\psit_N\|^2 \le V_0(t,\psit_N) &\le a_2\|\psit_N\|^2 \\
    \frac{\partial V_0}{\partial t} + \frac{\partial V_0}{\partial \psit_N}\big(-\gamma\xis_N(\xis_N)\T\psit_N\big) &\le -a_3\|\psit_N\|^2 \\
    \Big\|\frac{\partial V_0}{\partial \psit_N}\Big\| &\le a_4\|\psit_N\| \,.
\end{align}
\end{subequations} 

Let $r\0 := \|\xit(t_0)\|$.  Since $\Fbz$ is Hurwitz, there exist constants $k_2 \geq 1$ and $\lambda_F>0$
such that
\begin{align}\label{eq:xitbound}
\|\tilde{\xi}(t)\|
\leq k_2 e^{-\lambda_F(t-t_0)}
\|\tilde{\xi}(t_0)\|,
\qquad t\geq t_0.
\end{align}
Hence $\|\xit(t)\| \le k_2 r\0$ for all $t \ge t_0$. Let $P$ satisfy the Lyapunov equation $\Fbz\T P + P\Fbz = -I$, and let $k_3 := \|\xis\|_{\cLinf}$. Using the triangle inequality and sub-multiplicativity, we have
\begin{align}\label{eq:gbound}
    \|g(t, \psit_N, \psit_R, \xit)\| \le \gamma k_4 \|\xit\|\,\|\psit_N\| + k_5 \|\psit_R\| + k_6 \|\xit\| \,,
\end{align}
where
\begin{align*}
    k_4 &:= 2k_3 + k_2 r\0 \,, \\
    k_5 &:= \gamma\big(k_3 + k_2 r\0\big)^2 \,, \\
    k_6 &:= \gamma\big(k_3 + k_2 r\0\big)\|\tho + \psis\| \,.
\end{align*}
Define the constants
\begin{align*}
    k_7 &:= \frac{1}{2\mu}\Big(1 + \frac{a_4^2 k_5^2}{a_3}\Big)\,, \quad
    k_8 := 1 + \frac{a_4^2 k_6^2}{a_3} \\
    k_9 &:= \min\{a_1, k_7, k_8\lambda_{\min}(P)\} \,, \\
    k_{10} &:= \max\{a_2, k_7, k_8\lambda_{\max}(P)\} \,,\\
    \lambda &:= \min\Big\{\frac{a_3}{2}, 1\Big\} \,.
\end{align*}
Consider the Lyapunov function
\begin{align*}
    V(t,\psit_N,\psit_R,\xit) := V_0(t,\psit_N) + k_7\|\psit_R\|^2 + k_8\,\xit\T P\xit \,,
\end{align*}
which satisfies
\begin{align}\label{eq:Vsandwich}
    k_9\|(\psit_N,\psit_R,\xit)\|^2 \le V \le k_{10}\|(\psit_N,\psit_R,\xit)\|^2 \,.
\end{align}
Using \eqref{eq:converse} and \eqref{eq:gbound}, we obtain
\begin{align*}
    \dot{V} &\le -a_3\|\psit_N\|^2 + a_4\|\psit_N\|\,\|g\| - 2\mu k_7\|\psit_R\|^2 - k_8\|\xit\|^2 \\
    &\le -a_3\|\psit_N\|^2 + a_4\gamma k_4\|\xit\|\,\|\psit_N\|^2 + a_4 k_5\|\psit_N\|\,\|\psit_R\| \\
    &\quad + a_4 k_6\|\psit_N\|\,\|\xit\| - 2\mu k_7\|\psit_R\|^2 - k_8\|\xit\|^2 \,.
\end{align*}
Applying Young's inequality to the two cross terms,
\begin{align*}
    a_4 k_5\|\psit_N\|\,\|\psit_R\| &\le \frac{a_3}{4}\|\psit_N\|^2 + \frac{a_4^2 k_5^2}{a_3}\|\psit_R\|^2 \,, \\
    a_4 k_6\|\psit_N\|\,\|\xit\| &\le \frac{a_3}{4}\|\psit_N\|^2 + \frac{a_4^2 k_6^2}{a_3}\|\xit\|^2 \,.
\end{align*}
Using the definitions of $k_7$, $k_8$, and $\lambda$ we get
\begin{align*}
    \dot{V} &\le -\frac{a_3}{2}\|\psit_N\|^2 - \|\psit_R\|^2 - \|\xit\|^2 + a_4\gamma k_4\| \xit(t) \|\,\|\psit_N\|^2 \\
    &\le -\lambda\|(\psit_N,\psit_R,\xit)\|^2 + a_4\gamma k_4\|\xit(t)\|\,\|(\psit_N,\psit_R,\xit)\|^2 \\
    &\le \Big(-\frac{\lambda}{k_{10}} + \frac{a_4\gamma k_4}{k_9}\|\xit(t)\|\Big) V \,.
\end{align*}
By the Comparison Lemma \cite[Lemma~3.4]{KHALIL02},
\begin{align*}
    V(t) \le V(T_\rho) \exp\Big( -\frac{\lambda}{k_{10}}(t - T_\rho) + \frac{a_4\gamma k_4}{k_9} \int_{T_\rho}^{t} \|\xit(s)\|\, ds \Big) \,.
\end{align*}
By \eqref{eq:xitbound},
\begin{align*}
    \int_{T_\rho}^{\infty} \|\xit(s)\|\, ds \le \frac{k_2 r\0}{\lambda_F} \,,
\end{align*}
so
\begin{align}\label{eq:Vsecond}
    V(t) &\le V(T_\rho)\, e^{k_{11}} \, e^{-k_{12}(t - T_\rho)} \,, \quad \forall\, t \ge T_\rho \,, \\
    \nonumber
    k_{11} &:= \frac{a_4\gamma k_4 k_2 r\0}{k_9\lambda_F} \,, \quad k_{12} := \frac{\lambda}{k_{10}}\,.
\end{align}

\emph{Step 3: Combining the two intervals.} 
We now bound $V$ on $[t_0, T_\rho]$.  
By \eqref{eq:firstint} and \eqref{eq:xitbound}, for all $t \in [t_0, T_\rho]$,
\begin{align*}
    \|(\psit_N, \psit_R, \xit)(t)\|^2 = \|\psit(t)\|^2 + \|\xit(t)\|^2 \le k_1^2 + k_2^2 r\0^2 \,.
\end{align*}
By \eqref{eq:Vsandwich},
\begin{align}\label{eq:Vfirst}
    V(t) &\le k_{10}\big(k_1^2 + k_2^2 r\0^2\big) \,, 
    \quad \forall\, t \in [t_0, T_\rho] \,.
\end{align}
Define $k_{13} := k_{10}\big(k_1^2 + k_2^2 r\0^2\big) e^{k_{11} + k_{12}(T_\rho - t_0)}$. 
For $t \in [t_0, T_\rho]$ we have $e^{k_{11}} \ge 1$ and $e^{-k_{12}(t - T_\rho)} \ge 1$, so \eqref{eq:Vfirst} becomes: for all $t \in [t_0, T_\rho]$
\begin{align}\label{eq:Vall}
    V(t) &\le k_{13}\, e^{-k_{12}(t - t_0)} \,.
\end{align}
For $t \ge T_\rho$, using \eqref{eq:Vsecond} and \eqref{eq:Vfirst}, we have 
\begin{align*}
V(t) &\le k_{10}\big(k_1^2 + k_2^2 r\0^2\big) e^{k_{11}} e^{-k_{12}(t - T_\rho)} = k_{13} e^{-k_{12}(t - t_0)} \,.
\end{align*}
Hence \eqref{eq:Vall} holds for all $t \ge t_0$. Combining with \eqref{eq:Vsandwich}, we get $\psih(t) \to \psis$ exponentially. 

\bibliographystyle{plain}
\bibliography{SIDBib} 

\end{document}